\documentclass[psamsfonts]{amsart}
\usepackage{amssymb}
\usepackage[all]{xy}
\theoremstyle{plain}
\newtheorem{thm}{Theorem}[section]

\newtheorem{lemma}[thm]{Lemma}
\newtheorem{cor}[thm]{Corollary}
\newtheorem{question}[thm]{Question}

\newtheorem{hypo}[thm]{Hypothesis}

\theoremstyle{definition}
\newtheorem{dfn}[thm]{Definition}
\newtheorem{example}[thm]{Example}

\theoremstyle{remark}

\newcommand{\semid}{\unitlength.47cm
 \begin{picture}(.7,.6)
  \put(0,.05){$\times$}
  \put(.47,.04){\line(0,1){.39}}
 \end{picture}}

\begin{document}

\title{Deformation rings which are not local complete intersections}

\date{\today}
\author{Frauke M. Bleher}
\address{F.B.: Department of Mathematics\\University of Iowa\\
Iowa City, IA 52242-1419}
\email{fbleher@math.uiowa.edu}
\author{Ted Chinburg}
\address{T.C.: Department of Mathematics\\University of
Pennsylvania\\Philadelphia, PA
19104-6395}
\email{ted@math.upenn.edu}
\author{Bart de Smit}
\address{B.deS.:  Mathematisch Instituut\\University of Leiden\\P.O. Box 9512
\\ 2300 RA Leiden\\
The Netherlands}
\email{desmit@math.leidenuniv.nl}

\thanks{The first author was supported in part by  
NSF Grant  DMS0651332.
The second author was supported in part by  NSF Grant DMS0801030.
The third author was funded in part by the European Commission under contract
MRTN-CT-2006-035495.}

\begin{abstract}
We study the inverse problem for the versal deformation rings $R(\Gamma,V)$ of finite dimensional representations $V$
of a finite group $\Gamma$ over a field $k$ of positive characteristic $p$.  This problem is to determine which complete local commutative Noetherian
rings with residue field $k$ can arise up to isomorphism as such $R(\Gamma,V)$.  We show that for all integers $n \ge 1$
and all complete local commutative Noetherian rings $\mathcal{W}$ with residue field $k$, 
the ring $\mathcal{W}[[t]]/(p^n t,t^2)$ arises in this way.
This ring is not a local complete intersection if $p^n\mathcal{W}\neq\{0\}$, so we obtain an answer to a question of M. Flach
in all characteristics.
\end{abstract}

\maketitle

%%%%%%%%%%%%%%%%%%%%%%%%%%%%%%%%%%%%%%%%%%%%%%%%%%%%%%%%%%%%%%%%%%%%%%%%%%%
%% Introduction
%%%%%%%%%%%%%%%%%%%%%%%%%%%%%%%%%%%%%%%%%%%%%%%%%%%%%%%%%%%%%%%%%%%%%%%%%%%

\section{Introduction}
\label{s:intro}

Suppose $\Gamma$ is a profinite group and that $V$ is a continuous finite dimensional representation
of $\Gamma$ over a field $k$ of characteristic $p > 0$.  Let $\mathcal{W}$ be a complete local 
commutative Noetherian ring with residue field $k$.  In \S \ref{s:basechange} we recall the definition 
of a deformation of $V$ over a complete local commutative Noetherian $\mathcal{W}$-algebra with 
residue field $k$.  It follows from work of Mazur and Schlessinger \cite{maz1,Sch} 
that $V$ has a Noetherian versal deformation ring $R_{\mathcal{W}}(\Gamma,V)$
if the $p$-Frattini quotient of every open subgroup of $\Gamma$ is finite.   Without assuming
this condition, de Smit and Lenstra proved in \cite{desmitlenstra} that $V$ has a universal deformation 
ring $R_{\mathcal{W}}(\Gamma,V)$ if $\mathrm{End}_{k\Gamma}(V) = k$.  The ring 
$R_{\mathcal{W}}(\Gamma,V)$ is a pro-Artinian $\mathcal{W}$-algebra, but it need not be Noetherian.   
In this paper we consider the following inverse problem:

\begin{question}
\label{q:whatsup}
Which complete local commutative Noetherian $\mathcal{W}$-algebras $R$ with residue field $k$ are 
isomorphic to $R_{\mathcal{W}}(\Gamma,V)$ for some $\Gamma$ and $V$ as above? 
\end{question}

It is important to emphasize that in this question, $\Gamma$ and $V$ are not fixed. Thus for a given $R$,
one would like to construct both a profinite group $\Gamma$ and a continuous finite dimensional 
representation $V$ of $\Gamma$ over  $k$ for which $R_{\mathcal{W}}(\Gamma,V)$ is isomorphic to 
$R$.
We will be most interested in the case of finite groups $\Gamma$ in this paper, for which $R_{\mathcal{W}}(\Gamma,V)$ is always Noetherian.

Our main result is:

\begin{thm} 
\label{thm:vaguemain} For all fields $k$ and rings $\mathcal{W}$ as above, and for all $n \ge 1$,  there is a 
 representation $V$ of a finite group $\Gamma$ over $k$ having a universal deformation ring $R_{\mathcal{W}}(\Gamma,V)$ which is isomorphic to $\mathcal{W}[[t]]/(p^n t,t^2)$.  This ring is not a local complete 
intersection if $p^n \mathcal{W} \ne \{0\}$.
\end{thm}

Recall (see \cite[\S 19.3]{ega44}) that a commutative local Noetherian ring  $R$ is a local complete intersection  if there is a regular complete local commutative Noetherian ring $S$ and a regular sequence $x_1,\ldots,x_n\in S$ such that the completion $\hat{R}$ is isomorphic to 
$S/(x_1,\ldots, x_n)$.  
The problem of constructing  representations having universal deformation rings which 
are not local complete intersections was first posed
by M. Flach   \cite{flach}.  The first example of a representation of
this kind was found by Bleher and Chinburg with $k = \mathbb{Z}/2$;  see \cite{lcicomptes} 
and \cite{lcann}.
A more elementary argument proving the same result for $k = \mathbb{Z}/2$ was given
by Byszewski in \cite{JB}. 

Before outlining the proof of Theorem  \ref{thm:vaguemain} we discuss some other rings for which Question \ref{q:whatsup}
has been shown to have an affirmative answer.  We will suppose in this discussion that $k$ is a perfect 
field,  and we let $\mathcal{W}$ be the ring $W(k)$ of infinite Witt vectors over $k$.

Work of Mazur concerning $V$ of dimension
$1$ over $k$ shows that  $W(k)[[H]]$ is a versal deformation ring if $k$ has characteristic 
$p > 0$ and $H$ is a topologically finitely generated abelian pro-$p$-group. 
 
In \cite{BC}, Bleher and Chinburg
considered $V$ which belong to blocks with cyclic defect groups of
the group ring $k\Gamma$ of a finite group $\Gamma$ over an algebraically closed field $k$
of characteristic $p$. They showed that $W(k)$ and $W(k)/W(k)p^d$ are  universal deformation
rings for all integers $d \ge 1$.  Their results also show that if $D$ is a finite cyclic $p$-group,
and $E$ is a finite group of automorphisms of $D$ of order dividing $p-1$, then
the ring $(W(k)D)^E/W(k)s$ is a versal deformation ring, where $(W(k)D)^E$ is the ring
of $E$-invariants in the group ring $W(k)D$, $s = 0$ if $E$ is trivial and $s = \sum_{d \in D} d$
otherwise.  
 
In \cite{FB1,FB2}, Bleher considered $V$ which belong to certain blocks with dihedral defect groups,
respectively with generalized quaternion defect groups,
of the group ring $k\Gamma$ of a finite group $\Gamma$ over an algebraically closed field $k$
of characteristic $2$. She showed that $W(k)[[t]]/(p_d(t)(t-2),2\,p_d(t))$ and
$W(k)[[t]]/(p_{d}(t))$ are universal deformation rings for all integers $d\ge 3$,
where $p_d(t)$ is the product of the minimal polynomials of $\zeta_{\ell}+\zeta_{\ell}^{-1}$ over $W(k)$
for $\ell\in\{2,\ldots, d-1\}$ 
and $\zeta_{\ell}$ is a primitive $2^\ell$-th root of unity.
 
As of this writing we do not know of a complete local commutative Noetherian ring $R$ with perfect 
residue field $k$ of positive characteristic which cannot be
realized as a versal deformation ring of the form $R_{W(k)}(\Gamma,V)$ for some profinite $\Gamma$
and some representation $V$ of $\Gamma$ over $k$.  
 
We now describe the sections of this paper.  
 
In \S \ref{s:basechange} we recall the notations of deformations and of versal and universal
deformation rings. We show in Theorem \ref{thm:basechange} that versal deformation
rings respect arbitrary base changes, generalizing a result of Faltings concerning base
changes from $W(k)$ to $W(k')$ when $k'$ is a finite extension of a perfect field $k$.  
This reduces the proof of  Theorem \ref{thm:vaguemain} 
to the case in which $k = \mathbb{Z}/p$ and $\mathcal{W}  = W(k) = \mathbb{Z}_p$.
  
In \S \ref{s:computit} we consider arbitrary perfect fields $k$ of characteristic $p$ and
we take $\mathcal{W} = W(k)$. 
In Theorem \ref{thm:genresult} we give 
a sufficient set of conditions on a representation $\tilde {V}$ of a finite
group $\Gamma$ over $k$ for the universal deformation ring $R_{W(k)}(\Gamma,\tilde{V})$
to be isomorphic to $R = W(k)[[t]]/(p^n t,t^2)$.  The proof that these conditions
are sufficient involves first showing that $R_{W(k)}(\Gamma,\tilde{V})$ is a
quotient of $W(k)[[t]]$ by proving that the dimension of the tangent space
of the deformation functor associated to $\tilde{V}$ is one.  We then construct
an explicit lift of $\tilde{V}$ over $R$ and show that this
cannot be lifted further to any small extension ring of $R$ which is a quotient
of $W(k)[[t]]$. 
 
In \S \ref{s:semdis} we show that the hypotheses of Theorem \ref{thm:genresult} 
are satisfied in certain cases when $\Gamma$ is isomorphic to an iterated semi-direct product 
$V' \semid ((\mathbb{Z}/\ell) \semid (\mathbb{Z}/q))$ of an abelian $p$-group
 $V'$ with cyclic groups of orders $\ell$ and $ q $ greater than $1$ which are prime to $p$.   The resulting examples are sufficient to complete the proof of
 Theorem \ref{thm:vaguemain}.
 
\medbreak
\noindent {\bf Acknowledgements:}
The authors would like to thank
M. Flach for correspondence about
his question.  The second author would also like to thank
the University of Leiden for its hospitality during the
spring of 2009.
This paper will appear elsewhere in final form.

%%%%%%%%%%%%%%%%%%%%%%%%%%%%%%%%%%%%%%%%%%%%%%%%%%%%%%%%%%%%%%%%%%%%%%%%%%%
%% Deformation rings of base changes
%%%%%%%%%%%%%%%%%%%%%%%%%%%%%%%%%%%%%%%%%%%%%%%%%%%%%%%%%%%%%%%%%%%%%%%%%%%

\section{Deformation rings of base changes}
\label{s:basechange}

Let $\Gamma$ be a profinite group, and 
let $k$ be a field of characteristic $p > 0$.  Let $\mathcal{W}$ be a complete local commutative
Noetherian ring with residue field $k$.  We denote by $\hat{\mathcal{C}}$ the category of all 
complete local commutative Noetherian $\mathcal{W}$-algebras with residue field $k$.  
Homomorphisms in $\hat{\mathcal{C}}$ are continuous $\mathcal{W}$-algebra homomorphisms 
which induce the identity map on $k$.  
Define $\mathcal{C}$ to be the  full subcategory of Artinian objects in $\hat{\mathcal{C}}$. 
For each ring $A$ in $\hat{\mathcal{C}}$, let $\mathfrak{m}_A$ be its maximal ideal and denote
the surjective morphism $A\to A/\mathfrak{m}_A=k$ in $\hat{\mathcal{C}}$ by $\pi_A$.
If $\alpha:A\to A'$ is a morphism in  $\hat{\mathcal{C}}$, we denote the induced morphism
$\mathrm{GL}_n(A)\to \mathrm{GL}_n(A')$ also by $\alpha$.

Let $n$ be a positive integer, and let $\rho:\Gamma\to \mathrm{GL}_n(k)$ be a
continuous homomorphism, where $\mathrm{GL}_n(k)$ has the discrete topology. 
By a lift of $\rho$ over a ring $A$ in $\hat{\mathcal{C}}$ we mean a continuous homomorphism
$\tau:\Gamma\to \mathrm{GL}_n(A)$ such that $\pi_A\circ\tau=\rho$. We say two lifts 
$\tau,\tau':\Gamma\to \mathrm{GL}_n(A)$ of $\rho$ over $A$ are strictly equivalent
if one can be brought into the other by conjugation by a matrix in the kernel of
$\pi_A: \mathrm{GL}_n(A)\to  \mathrm{GL}_n(k)$. We call a strict equivalence class of
lifts of $\rho$ over $A$ a deformation of $\rho$ over $A$ and define $\mathrm{Def}_\rho(A)$ to be the 
set of deformations $[\tau]$ of lifts $\tau$ of $\rho$ over $A$. We then have a functor 
$$\hat{H}_\rho:\hat{\mathcal{C}}\to\mathrm{Sets}$$
which sends a ring $A$ in $\hat{\mathcal{C}}$ to the set $\mathrm{Def}_\rho(A)$. Moreover,
if $\alpha:A\to A'$ is a morphism in $\hat{\mathcal{C}}$, then 
$\hat{H}_\rho(\alpha):\mathrm{Def}_\rho(A)\to \mathrm{Def}_\rho(A')$  sends
a deformation $[\tau]$ of $\rho$ over $A$ to the deformation $[\alpha\circ\tau]$
of $\rho$ over $A'$.

Instead of looking at continuous matrix representations of $\Gamma$, we can also look at
topological $\Gamma$-modules as follows.
Let $V=k^n$ be endowed with the continuous $\Gamma$-action given by composition of $\rho$ with
the natural action of $\mathrm{GL}_n(k)$ on $V$, i.e. $V$ is the $n$-dimensional topological 
$k\Gamma$-module corresponding to $\rho$. A lift of $V$ over a ring $A \in \hat{\mathcal{C}}$ is then
a pair $(M,\phi)$ consisting of a finitely generated free $A$-module $M$ on which $\Gamma$ acts 
continuously together with a $\Gamma$-isomorphism $\phi:k \otimes_A M \to V$ of (discrete) $k$-vector 
spaces.  We define $\mathrm{Def}_V(A)$ to be the set of isomorphism classes $[M,\phi]$ of lifts 
$(M,\phi)$ of $V$ over $A$.  We then have a functor 
$$\hat{F}_V:\hat{\mathcal{C}}\to\mathrm{Sets}$$
which sends a ring $A$ in $\hat{\mathcal{C}}$ to the set $\mathrm{Def}_V(A)$. Moreover,
if $\alpha:A\to A'$ is a morphism in $\hat{\mathcal{C}}$, then 
$\hat{F}_V(\alpha):\mathrm{Def}_V(A)\to \mathrm{Def}_V(A')$  sends
a deformation $[M,\phi]$ of $V$ over $A$ to the deformation $[A'\otimes_{A,\alpha}M,\phi_\alpha]$
of $V$ over $A'$, where $\phi_\alpha$
is the composition $k\otimes_{A'}(A'\otimes_{A,\alpha}M)\cong k\otimes_AM\xrightarrow{\phi} V$. 
The functors $\hat{F}_V$ and $\hat{H}_\rho$ are naturally isomorphic.

One says that a ring $R=R_{\mathcal{W}}(\Gamma,\rho)$ (resp. $R=R_{\mathcal{W}}(\Gamma,V)$) in $\hat{\mathcal{C}}$ is a versal 
deformation ring for $\rho$ (resp. for $V$) if there is a lift $\nu:\Gamma\to\mathrm{GL}_n(R)$ of $\rho$ 
over $R$ (resp. a lift $(U,\phi_U)$ of $V$ over $R$) such that the following conditions hold.  
For all rings $A$ in $\hat{\mathcal{C}}$, the map 
$$f_A: \mathrm{Hom}_{\hat{\mathcal{C}}}(R,A)\to \mathrm{Def}_\rho(A)
\quad\mbox{(resp. $f_A: \mathrm{Hom}_{\hat{\mathcal{C}}}(R,A)\to \mathrm{Def}_V(A)$)}$$
which sends a morphism $\alpha:R\to A$ in $\hat{\mathcal{C}}$ to the deformation
$\hat{H}_\rho(\alpha)([\nu])$ (resp. $\hat{F}_V(\alpha)([U,\phi_U])$) is surjective.
Moreover, if $k[\epsilon]$ is the ring of dual numbers with $\epsilon^2=0$,
then $f_{k[\epsilon]}$ is bijective. (Here the $\mathcal{W}$-algebra structure of $k[\epsilon]$ 
is such that the maximal ideal of $\mathcal{W}$ annihilates $k[\epsilon]$.) We call the deformation $[\nu]$ (resp. $[U,\phi_U]$)
a versal deformation of $\rho$ (resp. of $V$) over $R$. 
By Mazur \cite[Prop. 20.1]{maz2}, $\hat{H}_\rho$ (resp. $\hat{F}_V$) is continuous, which means that
we only need to check the surjectivity of $f_A$ for Artinian rings $A$ in $\mathcal{C}$.
The versal deformation ring $R=R_{\mathcal{W}}(\Gamma,\rho)$ (resp. $R=R_{\mathcal{W}}(\Gamma,V)$) is unique up to isomorphism if it exists.

If the map $f_A$ is bijective for all rings $A$ in $\hat{\mathcal{C}}$, then we say
$R=R_{\mathcal{W}}(\Gamma,\rho)$ (resp. $R=R_{\mathcal{W}}(\Gamma,V)$) is a universal deformation ring of $\rho$ (resp. of $V$) and 
$[\nu]$ (resp. $[U,\phi_U]$) is a universal deformation of $\rho$ (resp. of $V$) over $R$. 
This is equivalent to saying that $R$ represents the 
deformation functor  $\hat{H}_\rho$ (resp. $\hat{F}_V$) in the sense that $\hat{H}_\rho$ 
(resp. $\hat{F}_V$)  is naturally isomorphic to the Hom functor 
$\mathrm{Hom}_{\hat{\mathcal{C}}}(R,-)$. 

We will suppose from now on that $\Gamma$ satisfies the following $p$-finiteness condition
used by Mazur in \cite[\S 1.1]{maz1}:

\begin{hypo}
\label{hyp:finite} For every open subgroup $J$ of finite index in $\Gamma$, there
are only a finite number of continuous homomorphisms from $J$ to $\mathbb{Z}/p$.
\end{hypo}

It follows by \cite[\S1.2]{maz1} that for $\Gamma$ satisfying Hypothesis \ref{hyp:finite}, all finite 
dimensional continuous representations $V$ of 
$\Gamma$ over $k$ have a versal deformation ring.  It is shown in  \cite[Prop. 7.1]{desmitlenstra} that
if $\mathrm{End}_{k\Gamma}(V)=k$, then $V$ has a universal deformation ring.

The following result shows how versal deformation rings change when extending the
residue field $k$. For finite extensions of $k$, this was proved by Faltings (see \cite[Ch. 1]{Wiles}).

\begin{thm}
\label{thm:basechange}
Let $\Gamma$, $k$, $\mathcal{W}$ and $\rho$ be as above.
Let $k'$ be a field extension of $k$, and let $\mathcal{W}'$ be a complete local commutative
Noetherian ring with residue field $k'$ such that there is a local homomorphism $\mathcal{W}\to 
\mathcal{W}'$.
Let $\rho':\Gamma\to \mathrm{GL}_n(k')$ be the composition of $\rho$ with the injection
$\mathrm{GL}_n(k)\hookrightarrow \mathrm{GL}_n(k')$. 
Then the versal deformation ring $R_{\mathcal{W}'}(\Gamma,\rho')$ is the completion
$R'$ of $\Omega=\mathcal{W}'\otimes_{\mathcal{W}}R_{\mathcal{W}}(\Gamma,\rho)$
with respect to the unique maximal ideal $\mathfrak{m}_\Omega$ of $\Omega$.
\end{thm}

\begin{proof} 
Define $\hat{\mathcal{C}}'$ to be the category of all 
complete local commutative Noetherian $\mathcal{W}'$-algebras with residue field $k'$,
and define $\mathcal{C}'$ to be the full subcategory of Artinian objects in $\hat{\mathcal{C}}'$.
Let $R=R_{\mathcal{W}}(\Gamma,\rho)$ and let $\nu:\Gamma\to\mathrm{GL}_n(R)$ be a versal lift of 
$\rho$ over $R$.
Define $\Omega=\mathcal{W}'\otimes_{\mathcal{W}}R$ and let $R'$ be the completion
of $\Omega$ with respect to $\mathfrak{m}_\Omega$.
Then $\nu':\Gamma\to\mathrm{GL}_n(R')$ is a lift of $\rho'$ over $R'$ when we let 
$\nu'(g)=\left(1\otimes\nu(g)_{i,j}\right)_{1\le i,j\le n}$ for all $g\in \Gamma$.
We will first show that if $A' \in \mathrm{Ob}(\mathcal{C}')$ is Artinian and 
$\tau':\Gamma\to\mathrm{GL}_n(A')$ is a lift of $\rho'$ over $A'$, then there is a morphism 
$\alpha:R'  \to A'$ in $\hat{\mathcal{C}}'$ such that 
\begin{equation}
\label{eq:firststep}
[\tau']=[\alpha\circ\nu'].
\end{equation}
In the following, when there is an implicit homomorphism $\phi:Y \to X$ of modules, we will abbreviate
$X/\phi(Y)$ by $X/Y$. 
Tensoring the short exact sequence
$$0\to \mathfrak{m}_R \to R \to k \to 0$$
with $\mathcal{W}'$ over $\mathcal{W}$ gives an isomorphism
$$\frac{\Omega}{\mathcal{W}' \otimes_{\mathcal{W}} \mathfrak{m}_R } \cong 
\mathcal{W}' \otimes_{\mathcal{W}} k .$$
Thus we have isomorphisms
$$\frac{\Omega}{\mathfrak{m}_{\mathcal{W}'} +\mathcal{W}' \otimes_{\mathcal{W}} \mathfrak{m}_R} \cong \frac{\mathcal{W}' \otimes_{\mathcal{W}} k}{\mathfrak{m}_{\mathcal{W}'} \otimes_{\mathcal{W}} k} 
\cong k' \otimes_k k \cong k' .$$
This proves that the natural homomorphism
\begin{equation}
\label{eq:maxomega}
\mathfrak{m}_{\mathcal{W}'} +\mathcal{W}' \otimes_{\mathcal{W}} \mathfrak{m}_R\to \mathfrak{m}_{\Omega} 
\end{equation}
is surjective.
Since $A'$ is (discrete) Artinian, (\ref{eq:maxomega}) implies that
\begin{equation}
\label{eq:seemtoneed}
\mathrm{Hom}_{\hat{\mathcal{C}}'}(R',A')=
\mathrm{Hom}_{\mathrm{cont}}(\Omega,A')
\end{equation}
where $\mathrm{Hom}_{\mathrm{cont}}$ stands for the space of continuous $\mathcal{W}'$-algebra
homomorphisms which induce the identity map on the residue field $k'$.

The kernel $J$ of $\rho'$ is open of finite index in $\Gamma$.  Because of Hypothesis
\ref{hyp:finite}, there is a finite subset $S \subset J$ which projects to a set of generators for the
$p$-Frattini quotient of $J$.  Then $S$ projects to a set of topological generators for the maximal
pro-$p$ quotient $J' = J/J_0$ of $J$, where $J_0$ is normal in $\Gamma$ because $J$
is normal in $\Gamma$ and $J_0$ is characteristic in $J$.    
Let $S_0$ be a finite set of coset representatives
for $J$ in $\Gamma$.   Since 
$\tau'(J)$ is a pro-$p$ group, the set $\{\tau'(g): g \in S \cup S_0\}$ is a finite set of 
topological generators for the image of $\tau'$.  Because $\rho'$ and $\rho$ have the same 
image in $\mathrm{GL}_n(k)\subset\mathrm{GL}_n(k')$ and because $\tau'$ is a 
lift of $\rho'$, the following is true for each element $g$ of the set $S \cup S_0$.
There is a matrix $t(g) \in \mathrm{Mat}_n(\mathcal{W})$  such that all entries of
the matrix $\tau'(g) - t(g)$ lie in the maximal ideal $\mathfrak{m}_{A'}$ of $A'$. 
Let $T$ be the finite set of all matrix entries of $\tau'(g) - t(g)$ as $g$  ranges over $S \cup S_0$,
so that $T \subset \mathfrak{m}_{A'}$.   Now every finite word in the elements of $S \cup S_0$
is sent by $\tau'$ to a matrix with entries in the $\mathcal{W}$-subalgebra of $A'$ generated
by $T$.  Since the elements of $T$ are in $\mathfrak{m}_{A'}$, there is a 
continuous homomorphism $f:\mathcal{W}[[x_1,\ldots,x_m]] \to A'$
with $m = \# T$ and $\{f(x_i)\}_{i = 1}^m = T$.   Let $B$ be the image of $f$ in $A'$.
Since $\{\tau'(g): g \in S \cup S_0\}$ is a finite set of 
topological generators for the image of $\tau'$ and since the matrix entries of $\tau'(g)$ for
$g\in S\cup S_0$ lie in $B$, it follows that the image of $\tau'$ lies in $\mathrm{GL}_n(\overline{B})$
when $\overline{B}$ is the closure of $B$ in $A'$.  However, $\overline{B} = B$ since $A'$ has
the discrete topology.
Because $\mathcal{W}[[x_1,\ldots,x_m]]$ is  a complete local commutative Noetherian ring with 
residue field equal to the residue field $k$ of $\mathcal{W}$, we have that 
$B=f(\mathcal{W}[[x_1,\ldots,x_m]])$ is a complete local commutative Noetherian ring with 
residue field $k$.  
Thus $\tau'$ defines a lift of $\rho$
over $B$ since $\rho'$ and $\rho$ have the same image and the image
of $\tau'$ lies in $\mathrm{GL}_n(B)$.   Because   $\nu:\Gamma\to\mathrm{GL}_n(R)$ 
is a versal lift of $\rho$ over the versal deformation ring $R = R_{\mathcal{W}}(\Gamma,\rho)$ of $\rho$,
there is a  morphism 
$\beta:R  \to B$ in $\hat{\mathcal{C}}$ such that $\tau':\Gamma \to \mathrm{GL}_n(B)$
is conjugate to $\beta \circ \nu$ by a matrix in the kernel of
$\pi_B:\mathrm{GL}_n(B) \to \mathrm{GL}_n(B/\mathfrak{m}_B) = \mathrm{GL}_n(k)$.
Let $\beta':R \to A'$ be the composition of $\beta$ with the inclusion $B \subset A'$.
Define $\alpha: R'  \to A'$
to be the morphism in $\hat{\mathcal{C}}'$ corresponding by (\ref{eq:seemtoneed})
to the continuous $\mathcal{W}'$-algebra homomorphism 
$\Omega=\mathcal{W}' \otimes_{\mathcal{W}} R  \to A'$ which sends $w'\otimes r$ to $w'\cdot\beta'(r)$
for all $w'\in\mathcal{W}'$ and $r\in R$.
Then $\alpha$ has the required property (\ref{eq:firststep}).

We must still show that when $k'[\epsilon]$ is the ring of dual numbers over
$k'$, then $\mathrm{Hom}_{\hat{\mathcal{C}}'}(R',k'[\epsilon])$
is canonically identified with the set $\mathrm{Def}_{\rho'}(k'[\epsilon])$ of deformations of $\rho'$ over
$k'[\epsilon]$.  
Since $\mathrm{Ad}(\rho') = k' \otimes_k \mathrm{Ad}(\rho)$, we have from \cite[Prop. 21.1]{maz2} that there are canonical identifications 
\begin{equation}
\label{eq:H1com}
\mathrm{Def}_{\rho'}(k'[\epsilon]) = H^1(\Gamma,\mathrm{Ad}(\rho')) = 
k' \otimes_k H^1(\Gamma,\mathrm{Ad}(\rho))=k'\otimes_k \mathrm{Def}_{\rho}(k[\epsilon]).
\end{equation}
By (\ref{eq:seemtoneed}), we wish to show $\mathrm{Hom}_{\mathrm{cont}}(\Omega,k'[\epsilon])$ is 
identified with $\mathrm{Def}_{\rho'}(k'[\epsilon])$.  
As before, we will write $X/Y$ for a quotient module $X/\phi(Y)$ when there is an implicit homomorphism 
$\phi:Y \to X$.  Let 
\begin{equation}
\label{eq:quot}
T(\mathcal{W}',\Omega) = \frac{\mathfrak{m}_{\Omega}}{\mathfrak{m}_{\Omega}^2 + \Omega \cdot \mathfrak{m}_{\mathcal{W}'}}
\quad\mbox{and}\quad
T(\mathcal{W},R) = \frac{\mathfrak{m}_R}{\mathfrak{m}_R^2 + R \cdot \mathfrak{m}_{\mathcal{W}}}
\end{equation}
so that we have natural isomorphisms 
\begin{equation}
\label{eq:natural}
\mathrm{Hom}_{\mathrm{cont}}(\Omega,k'[\epsilon]) \cong \mathrm{Hom}_{k'}(T(\mathcal{W}',\Omega),k')\quad \mbox{and}\quad 
\mathrm{Hom}_{\hat{\mathcal{C}}}(R,k[\epsilon]) \cong \mathrm{Hom}_{k}(T(\mathcal{W},R),k).
\end{equation}
It follows from (\ref{eq:H1com}), (\ref{eq:quot}) and (\ref{eq:natural}) that to 
complete the proof, it will suffice to show that the natural homomorphism
\begin{equation}
\label{eq:mudef}
\mu:k' \otimes_k  T(\mathcal{W},R) \to T(\mathcal{W}',\Omega)
\end{equation}
is an isomorphism of $k'$-vector spaces.  
It follows from (\ref{eq:maxomega}) that the natural homomorphism
\begin{equation}
\label{eq:surjmust}
\mathcal{W}'\otimes_{\mathcal{W}} \mathfrak{m}_R \to 
\mathfrak{m}_{\Omega}/\mathfrak{m}_{\mathcal{W}'}
\end{equation}
is surjective, which implies that $\mu$ is surjective.
Hence to show  that $\mu$ is an isomorphism it will be enough to show
\begin{equation}
\label{eq:dims}
\mathrm{dim}_k T(\mathcal{W},R) = \mathrm{dim}_{k'} T(\mathcal{W}',\Omega).
\end{equation}
Suppose $\alpha \in \mathfrak{m}_\mathcal{W}$.  Define $\mathcal{W}_0 = \mathcal{W}/(\mathcal{W}\alpha)$
and $R_0 = R/(R\alpha)$.  Then $R_0$ is a complete local ring with maximal ideal $\mathfrak{m}_{R_0} = \mathfrak{m}_R/R \alpha$.
It follows that the natural projection
$T(\mathcal{W},R) \to T(\mathcal{W}_0,R_0)$
is an isomorphsm.  Similarly, 
if we define $\alpha' = \alpha \otimes 1  = 1 \otimes \alpha\in \mathcal{W}' \otimes_{\mathcal{W}} R=\Omega$,
then $T(\mathcal{W}',\Omega)$ does not change if we replace $\mathcal{W}'$ by $\mathcal{W}'_0 = \mathcal{W}'/(\mathcal{W}' \alpha)$ and $\Omega$ by $\Omega/(\Omega\alpha')$. Here
$$\Omega/(\Omega\alpha') \cong \mathcal{W}' \otimes_{\mathcal{W}} (R/R\alpha) \cong
\mathcal{W}'_0 \otimes_{\mathcal{W}_0} \mathcal{R}_0,$$
where the first isomorphism is from the right exactness of tensor products and the second
isomorphism is from the universal property of tensor products.
Thus to show (\ref{eq:dims}) we can  first divide all the rings involved by any ideal generated by an element of
$\mathfrak{m}_\mathcal{W}$.  Since $\mathfrak{m}_{\mathcal{W}}$ is finitely generated, we can thus reduce to the  case in which  $\mathcal{W} = k$.  We now divide $\mathcal{W}'$
and $\Omega$ further by ideals generated by generators for $\mathfrak{m}_{\mathcal{W}'}$ to be able to assume that
$\mathcal{W}' = k'$.  Thus we have reduced to the case in which $\mathfrak{m}_{\mathcal{W}} = \{0\} = \mathfrak{m}_{\mathcal{W}'}$,
$R$ is a complete local commutative $k$-algebra with residue field $k$ and 
$\Omega = k' \otimes_k R$.  
It follows that $\mathfrak{m}_{\Omega} = k' \otimes_k \mathfrak{m}_R$, and this identification
sends $k' \otimes_k \mathfrak{m}_{R}^2$ onto $\mathfrak{m}_{\Omega}^2$ inside 
$\mathfrak{m}_{\Omega}$.
Since we have a short exact sequence
$$0 \to k' \otimes_k \mathfrak{m}_{R}^2 \to k' \otimes_k \mathfrak{m}_R \to 
k' \otimes_k 
\left(\mathfrak{m}_R/\mathfrak{m}_R^2\right)
\to 0,$$
it follows that we have an isomorphism
$$T(k',\Omega) = \frac{\mathfrak{m}_{\Omega}}{\mathfrak{m}_{\Omega}^2} \cong k' \otimes_k 
\left(\frac{\mathfrak{m}_{R}}{\mathfrak{m}_{R}^2}\right)=  k' \otimes T(k,R).$$
This completes the proof of (\ref{eq:dims}) and of Theorem \ref{thm:basechange}.
\end{proof}

\begin{cor}
\label{cor:reduce} 
If for each prime $p$  there is a finite group $\Gamma$ and a finite dimensional
representation $T_0$ of $\Gamma$ over $\mathbb{F}_p$ such that
\begin{enumerate}
\item[(a)] $\mathrm{End}_{\mathbb{F}_p\Gamma}(T_0) = \mathbb{F}_p$, and 
\item[(b)] $R_{\mathbb{Z}_p}(\Gamma,T_0)$ is isomorphic to $\mathbb{Z}_p[[t]]/(p^n t,t^2)$
as an algebra over the ring  $W(\mathbb{F}_p) = \mathbb{Z}_p$
of $p$-adic integers,
\end{enumerate}
then Theorem $\ref{thm:vaguemain}$ holds.
\end{cor}

\begin{proof}  As in Theorem \ref{thm:vaguemain}, let $\mathcal{W}$ be a complete local
commutative Noetherian ring with residue field $k$ of characteristic $p$.  Consider 
the representation $T = k\otimes_{\mathbb{F}_p} T_0$ of $\Gamma$ over $k$.
Since  $$\mathrm{End}_{k\Gamma}(T) = 
k \otimes_{\mathbb{F}_p} \mathrm{End}_{\mathbb{F}_p\Gamma}(T_0) = k$$ 
an argument of Faltings (see \cite[\S 2.6]{Darmon} and \cite[\S7]{desmitlenstra}) shows
$T$ has a  universal deformation ring $R_{\mathcal{W}}(\Gamma,T)$ which is a complete local
commutative Noetherian  $\mathcal{W}$-algebra with residue field $k$.
By Theorem \ref{thm:basechange} and Hypothesis (b) of the Corollary we have isomorphisms 
$$R_{\mathcal{W}}(\Gamma,T) \cong \mathcal{W} \hat{\otimes}_{\mathbb{Z}_p} 
R_{\mathbb{Z}_p}(\Gamma,T_0) \cong \mathcal{W} \hat{\otimes}_{\mathbb{Z}_p} 
\mathbb{Z}_p[[t]]/(p^n t,t^2) \cong  \mathcal{W}[[t]]/(p^nt,t^2).$$  
Note that the Krull dimension of this ring is
\begin{equation}
\label{eq:krull}
\mathrm{dim}\, \mathcal{W}[[t]]/(p^nt,t^2) = \mathrm{dim}\, \mathcal{W}[[t]] - 1 = 
\mathrm{dim}\,\mathcal{W}.
\end{equation}
Suppose now $p^n\mathcal{W}\neq\{0\}$, and assume that $\mathcal{W}[[t]]/(p^nt,t^2)$ is a 
local complete intersection. Then $\mathcal{W}[[t]]/(p^nt,t^2)\cong S/J$ where $S$ is a regular
complete local commutative Noetherian ring and $J$ is an ideal of $S$ which is
generated by a regular sequence. Since $\mathcal{W}$ is a quotient of 
$\mathcal{W}[[t]]/(p^nt,t^2)$, $\mathcal{W}$ is also a quotient of $S$, say $\mathcal{W}\cong S/I$,
where $I$ is contained in the maximal ideal $\mathfrak{m}_S$.   The ring $S'=S[[t]]$ is a regular
complete local commutative Noetherian ring with maximal ideal $\mathfrak{m}_{S'}$. 
Let $I'$ be the ideal of $S'$ generated by $I$, $p^nt$ and $t^2$,
so $\mathcal{W}[[t]]/(p^nt,t^2)\cong S'/I'$.  Since $\mathcal{W}\cong S/I$ and we assumed $\mathcal{W}[[t]]/(p^nt,t^2)$ 
to be a local complete intersection, this implies by \cite[Thm. 21.1]{matsumura} and by (\ref{eq:krull})
that
\begin{eqnarray}
\mathrm{dim}_k(I'/\mathfrak{m}_{S'}I') &=& \mathrm{dim}\, S' -\mathrm{dim}\, (S'/I')\nonumber\\
&=& \mathrm{dim}\,S +1 - \mathrm{dim}\, (S/I)\nonumber\\
\label{eq:equal}
&\le &\mathrm{dim}_k(I/\mathfrak{m}_SI)+1.
\end{eqnarray}
Using power series expansions, we see that 
$I'=I + t \,(p^nS+I) + t^2\,S[[t]]$, $\mathfrak{m}_{S'} = \mathfrak{m}_S + tS[[t]]$ and $\mathfrak{m}_{S'}I'=\mathfrak{m}_SI + t\,(p^n\mathfrak{m}_S+I) + t^2\,\mathfrak{m}_S + t^3\,S[[t]]$.
Hence 
$$\frac{I'}{\mathfrak{m}_{S'} I'} \cong \frac{I}{\mathfrak{m}_S I} \oplus  \frac{p^nS+I}{p^n\mathfrak{m}_S+I} \oplus \frac{S}{\mathfrak{m}_S}.$$
Since $\mathcal{W}=S/I$ and $p^n\mathcal{W}\neq \{0\}$, it follows that 
$(p^nS+I)/(p^n\mathfrak{m}_S+I)\cong (p^n\,\mathcal{W})/(p^n\,\mathfrak{m}_{\mathcal{W}})\cong k$.
We obtain $\mathrm{dim}_k(I'/\mathfrak{m}_{S'}I') =\mathrm{dim}_k(I/\mathfrak{m}_SI)+2$,
which contradicts (\ref{eq:equal}). Thus $\mathcal{W}[[t]]/(p^nt,t^2)$ is not a local 
complete intersection if $p^n\mathcal{W}\neq\{0\}$.

\end{proof}

%%%%%%%%%%%%%%%%%%%%%%%%%%%%%%%%%%%%%%%%%%%%%%%%%%%%%%%%%%%%%%%%%%%%%%%%%%%
%% Computing deformation rings 
%%%%%%%%%%%%%%%%%%%%%%%%%%%%%%%%%%%%%%%%%%%%%%%%%%%%%%%%%%%%%%%%%%%%%%%%%%%

\section{Computing deformation rings} 
\label{s:computit} 
  
Throughout this section we make the following assumptions.

\begin{hypo}
\label{hyp:setup}
Let  $k$ be an arbitrary perfect field of characteristic $p > 0$
and let  $\mathrm{W}$ be the ring $W(k)$ of infinite Witt vectors over $k$. 
Let $G$ be a finite group and suppose $n \ge 1$ is an integer.  Define $A = \mathrm{W}/(\mathrm{W}p^n)$. 
Suppose $V$ is a projective $kG$-module
for which $\mathrm{dim}_k(V)$ is finite and $\mathrm{End}_{kG}(V) = k$.  Let $\hat{V}$ be a projective
$AG$-module such that $k \otimes_{A} \hat{V}$ is isomorphic to $V$
as a $kG$-module.  Let $M$ be the
free $A$-module $\mathrm{Hom}_A(\hat{V},\hat{V})$, so that $M$ is a projective $AG$-module.  Define
$$M_0 = k \otimes_A M = \mathrm{Hom}_k(V,V).$$
If $L$ is an $AG$-module,
we will also view $L$ as an $(\mathbb{Z}/p^n)G$-module via restriction of operators from 
 $AG$ to $(\mathbb{Z}/p^n)G$.  
 \end{hypo}
 
\begin{thm}
\label{thm:genresult}    
Assume Hypothesis $\ref{hyp:setup}$. Suppose $V'$ is a $(\mathbb{Z}/p^n)G$-module which is a free, 
finitely generated $\mathbb{Z}/p^n$-module with the following properties:
\begin{enumerate}
\item[(i)] The group $\mathrm{Hom}_{(\mathbb{Z}/p^n)G}(V',M)$ is a free rank one $A$-module with
respect to the $A$-module structure coming from the multiplication action of $A$ on $M$.
\item[(ii)]  There is an injective homomorphism $\psi:V' \to M$ in
$\mathrm{Hom}_{(\mathbb{Z}/p^n)G}(V',M)$.
\item[(iii)]  There are elements $\tau,\lambda \in \psi(V') \subset M = \mathrm{Hom}_A(\hat{V},\hat{V})$ and an element
$v \in \hat{V}$ such that $\tau(\lambda(v)) - \lambda(\tau(v)) \not \in p\hat{V}$.  
\end{enumerate}
Let $K=V'$ as abelian $p$-groups and let $\delta:G \to \mathrm{Aut}(K)$ be the group 
homomorphism given by
the $G$-action on the $(\mathbb{Z}/p^n)G$-module $V'$.
Define $\Gamma$ to be the semi-direct product $K \semid_{\delta} \,G$.
Let $\tilde{V}$ be the $\Gamma$-module which results by inflating the $G$-module $V$
via the natural surjection $\pi:\Gamma \to G$.  Then under these hypotheses,
the universal deformation ring $R_{\mathrm{W}}(\Gamma,\tilde{V})$ is well defined and is isomorphic
to $\mathrm{W}[[t]]/(p^nt,t^2)$.  
\end{thm}

\begin{example}
\label{eq:easyex}
Suppose $p = 2$,  that $k = \mathbb{F}_2$ and that 
$\mathrm{W}=W(\mathbb{F}_2)=\mathbb{Z}_2$. Let $G$ be the
symmetric group on $3$ elements.  The two dimensional irreducible representation
$V$ of $G$ over $k$ is projective.  Let $\hat{V}$ be a projective $AG$-module such
that $\hat{V}/p\hat{V}$ is isomorphic to $V$ as a $kG$-module, and let $V' = \hat{V}=K$.
Then the hypotheses of Theorem \ref{thm:genresult} are satisfied,
as may be seen in the following
way.  As a $kG$-module, $M_0 = M/pM$ is isomorphic to $V \oplus k[G/C]$ when $C$ is
the index two subgroup of $G$. We obtain from this an injection $\psi$ as in 
condition (ii) of Theorem \ref{thm:genresult}.  The image of $\psi$ mod $pM$ is a
two-dimensional subspace of $M_0$ which does not contain the identity element
of $M_0 = \mathrm{Hom}_k(V,V)$. Thus if condition (iii) failed, the identity element
together with the image of $\psi$ mod $pM$ would generate a commutative subalgebra
of $M_0$ of dimension at least $3$ over $k$, and no such subalgebra exists.  We thus obtain
from  Theorem \ref{thm:genresult} an extension $\Gamma$ of $G$ by $K$ such that $R_{\mathrm{W}}(\Gamma,\tilde{V})$ is  isomorphic
to $\mathrm{W}[[t]]/(p^nt,t^2)$.  When $n = 1$, 
the group $\Gamma$ is isomorphic to the symmetric group $S_4$ on four letters. 
The fact that $R_{\mathrm{W}}(\Gamma,\tilde{V})$  is isomorphic
to $\mathrm{W}[[t]]/(pt,t^2)$ in this case was shown in \cite{lcicomptes} 
and \cite{lcann};  see also  \cite{JB}.
\end{example}

We now return to Theorem \ref{thm:genresult}.  For any $G$-module $L$, we denote by
$\tilde{L}$ the $\Gamma$-module which results by inflating $L$
via the natural surjection $\pi:\Gamma \to G$. 

\begin{lemma}
\label{lem:cohomcor}
One has $\mathrm{dim}_k(H^1(\Gamma,\tilde{M_0})) = 1$.
The tangent space of the universal deformation ring $R_{\mathrm{W}}(\Gamma,\tilde{V})$ of $\tilde{V}$ has dimension
$1$.  The  ring $R_{\mathrm{W}}(\Gamma,\tilde{V})$
is a quotient of $\mathrm{W}[[t]]$.
\end{lemma}

\begin{proof} By construction there is an exact sequence of groups
$$1 \to K \to \Gamma \xrightarrow{\pi} G \to 1$$
in which $G$ acts on $K$ via $\delta$. Since $M_0$ is a projective $kG$-module, we have $H^i(G,H^0(K,\tilde {M_0})) = H^i(G,M_0) = 0$
if $i > 0$.  Therefore the Hochschild-Serre spectral sequence for $H^1(\Gamma,\tilde{M_0})$ degenerates to give 
\begin{equation}
\label{eq:stupid}
H^1(\Gamma,\tilde{M_0}) =H^0(G,H^1(K,\tilde {M_0})) = H^0(G,\mathrm{Hom}(K,\tilde {M_0})) = 
\mathrm{Hom}(K,M_0)^G.
\end{equation}
Since $M_0$ is an elementary abelian $p$-group, we have from condition (i) 
of Theorem \ref{thm:genresult}  that 
\begin{eqnarray}
\label{eq:stupider}
 \mathrm{Hom}(K,M_0)^G &=& \mathrm{Hom}(V',M_0)^G = \mathrm{Hom}(V'/pV',M_0)^G \nonumber\\
 &=& \left ((\mathbb{Z}/p) \otimes_{\mathbb{Z}/p^n} \mathrm{Hom}_{(\mathbb{Z}/p^n)}(V',M)\right )^G\nonumber\\
 &=&   (\mathbb{Z}/p) \otimes_{\mathbb{Z}/p^n} \mathrm{Hom}_{(\mathbb{Z}/p^n)G}(V',M) \cong  (\mathbb{Z}/p) \otimes_{\mathbb{Z}/p^n} A = k.
 \end{eqnarray}
On putting together
(\ref{eq:stupid}) and (\ref{eq:stupider}), we conclude
from  \cite[Prop. 21.1]{maz2} that there
is a natural isomorphism
\begin{equation}
\label{eq:tanspace}
t_{\tilde V} =_{\mathrm{def}} \mathrm{Hom}_k\left (\frac{\mathfrak{m}}{\mathfrak{m}^2 + p R_{\mathrm{W}}(\Gamma, \tilde{V})} ,k\right ) \to H^1(\Gamma,\tilde{M_0}) = k
\end{equation}
where  $t_{\tilde V}$ is the tangent space of the deformation functor of $\tilde{V}$, $\mathfrak{m}$ is the maximal ideal of the universal deformation ring $R_{\mathrm{W}}(\Gamma, \tilde{V})$
and $R_{\mathrm{W}}(\Gamma,\tilde{V})$ is a complete local $\mathrm{W}$-algebra.  This implies 
$$\mathrm{dim}_k \left (\frac{\mathfrak{m}}{\mathfrak{m}^2 + p R_{\mathrm{W}}(\Gamma, \tilde{V})} \right ) = 1$$
so there is a continuous surjection of $\mathrm{W}$-algebras $\mathrm{W}[[t]] \to R_{\mathrm{W}}(\Gamma,\tilde{V})$.
\end{proof}

We now construct an explicit lift of a matrix representation for the $k\Gamma$-module $\tilde{V}$.

Let $q = \mathrm{dim}_k(V)$. By assumption, $V$ is 
projective as a $kG$-module which is the reduction mod $p$ of a projective $AG$-module $\hat{V}$.   Thus there
is a matrix representation 
$\rho_{\mathrm{W}}:G \to \mathrm{GL}_q(\mathrm{W})$ whose reduction mod $p^n \mathrm{W}$ is a matrix representation
$\hat{\rho}:G \to \mathrm{GL}_q(A)$ for $\hat{V}$, and whose reduction mod $p\mathrm{W}$ is a matrix representation 
$\overline{\rho}:G \to \mathrm{GL}_q(k)$  for $V$.

Let $R = \mathrm{W}[[t]]/(p^n t,t^2)$.  We have an exact sequence of multiplicative groups
\begin{equation}
\label{eq:Rex}
1 \to (1 + t\mathrm{Mat}_q(R))^* \to \mathrm{GL}_q(R) \to \mathrm{GL}_q(\mathrm{W}) \to 1
\end{equation}
resulting from the natural isomorphism $R/tR = \mathrm{W}$.  The  isomorphism
$tR \to A = \mathrm{W}/p^n \mathrm{W}$ defined by $t w \to w \ \mathrm{mod} \ p^n\mathrm{W}$ for $w \in \mathrm{W} \subset R$ gives rise to isomorphisms of  groups 
\begin{equation}
\label{eq:kerfix}
(1 + t\mathrm{Mat}_q(R))^* \cong \mathrm{Mat}_q(A)^+ \cong M \cong \mathrm{Hom}_A(\hat{V},\hat{V})
\end{equation}
where $\mathrm{Mat}_q(A)^+$ is the additive group of $\mathrm{Mat}_q(A)$.
The conjugation action of $\rho_\mathrm{W}(G) \subset \mathrm{GL}_q(\mathrm{W})$ on $ (1 + t\mathrm{Mat}_q(R))^*$ which
results from (\ref{eq:Rex}) factors through the homomorphism $\rho_\mathrm{W}(G) \to \hat{\rho}(G) \subset \mathrm{GL}_q(A) = \mathrm{Aut}_A(\hat{V})$.
This action coincides with the action of $G$ on $M = \mathrm{Hom}_A(\hat{V},\hat{V})$ in (\ref{eq:kerfix}) coming from
the action of $G$ on $\hat{V}$ via $\hat{\rho}:G \to \mathrm{GL}_q(A)$.  

We have supposed in Theorem \ref{thm:genresult} that there is an injective $(\mathbb{Z}/p^n)G$-module
homomorphism $\psi:K=V' \to M$ which is unique
up to multiplication by an element of $A^*$ acting on $M$. 
We conclude that
there is a group homomorphism $\rho_R$ which makes the following diagram commutative:
\begin{equation}
\label{eq:lifter}
\xymatrix {
1\ar[r]&K \ar[d]^{\psi} \ar[r]& \Gamma \ar[d]^{\rho_R} \ar[r]^{\pi}&G\ar[d]^{\rho_\mathrm{W}}\ar[r]&1\\
1\ar[r]&M \ar[r]^{\nu}&\mathrm{GL}_q(R)\ar[r] &\mathrm{GL}_q(\mathrm{W})\ar[r]&1
}
\end{equation}
Here $\nu$ in the bottom row results from  (\ref{eq:Rex}) and (\ref{eq:kerfix}).
Since $\rho_\mathrm{W}$ is a lift of the matrix representation $\overline{\rho}:G \to \mathrm{GL}_q(k) = \mathrm{Aut}_k(V)$
over $\mathrm{W}$, we find that $\rho_R$ is a lift of $\overline{\rho}\circ \pi$ over $R$. Here
$ \overline{\rho}\circ \pi$ is a matrix representation for the $q$-dimensional 
representation $\tilde {V}$ of $\Gamma$ over $k$ which is inflated from the representation $V$
of $G$.

\begin{lemma}
\label{lem:wazzup}
Let $\gamma: R_{\mathrm{W}}(\Gamma,\tilde{V}) \to R = \mathrm{W}[[t]]/(p^nt,t^2)$ be the unique continuous $\mathrm{W}$-algebra homomorphism corresponding to the isomorphism class of the lift $\rho_R$ of $\tilde{V}$.  Then $\gamma$ is surjective.    There is a $\mathrm{W}$-algebra surjection $\mu:\mathrm{W}[[t]] \to R_{\mathrm{W}}(\Gamma,\tilde{V})$ whose composition with $\gamma$
is the natural surjection $\mathrm{W}[[t]] \to R = \mathrm{W}[[t]]/(p^nt,t^2)$.  The kernel of $\mu$ is an ideal of $\mathrm{W}[[t]]$ contained
in $(p^n t,t^2)$.  The homomorphism $\gamma$ is not an isomorphism if and only if there is a $\mathrm{W}[[t]]$-ideal $I\subset (p^n t,t^2)$ with the following properties:
\begin{enumerate}
\item[(a)] $(p^n t,t^2)/I$ is isomorphic to $k$
\item[(b)] There is a lift of $\rho_R$ to a group homomorphism $\rho_I:\Gamma \to \mathrm{GL}_q(\mathrm{W}[[t]]/I)$.
\end{enumerate}
\end{lemma}

\begin{proof} The ring $k[\epsilon]$ of dual numbers over $k$ is isomorphic to $R/pR = k[[t]]/(t^2)$, and 
$\gamma$ is surjective if and only if it induces a surjection 
\begin{equation}
\label{eq:tan}
\overline{\gamma}:\quad
\frac{R_{\mathrm{W}}(\Gamma,\tilde{V})}{\mathfrak{m}^2 + pR_{\mathrm{W}}(\Gamma,\tilde{V})} 
\longrightarrow \frac{R}{\mathfrak{m}^2_R + pR} = \frac{R}{pR}
\end{equation}
where $\mathfrak{m}$ is the maximal ideal of $R_{\mathrm{W}}(\Gamma,\tilde{V})$.  
If $\overline{\gamma}$ is not surjective, its image is $k$. Thus
to prove that $\gamma$ is surjective, it will suffice to show that the composition 
$\rho_{R/pR}$ of $\rho_R$
with the natural surjection $\mathrm{GL}_q(R) \to \mathrm{GL}_q(R/pR) = \mathrm{GL}_q(k[\epsilon])$ is not a matrix
representation of the trivial lift of $\tilde{V}$ over $k[\epsilon]$. However, the kernel of the action of $\Gamma$ on this trivial lift is $K \subset \Gamma$,
while $\rho_{R/pR}$ is not trivial on $K$, so $\gamma$ must be surjective. 

The tangent space of the deformation functor of $\tilde{V}$ is one dimensional by Lemma
\ref{lem:cohomcor}, so (\ref{eq:tan}) is in fact an isomorphism.
Let $r$ be any element of $R_{\mathrm{W}}(\Gamma,\tilde{V})$ such that $\gamma(r)$ is the class of
$t$ in $R = \mathrm{W}[[t]]/(p^n t,t^2)$.  We then have a unique continuous $\mathrm{W}$-algebra homomorphism
$\mu:\mathrm{W}[[t]] \to R_{\mathrm{W}}(\Gamma,\tilde{V})$ which maps $t$ to $r$.  
Since $(\gamma\circ\mu)(t)$ is the class of $t$ in $R$, we se that $\gamma\circ\mu$ is surjective.
So because $\overline{\gamma}$ is an isomorphism, Nakayama's lemma implies
that $\mu:\mathrm{W}[[t]] \to R_{\mathrm{W}}(\Gamma,\tilde{V})$ is surjective.  The rest of
the Lemma \ref{lem:wazzup} is clear, since the kernel of $\mu$ is smaller than $(p^n t,t^2)$ if and only if there is an ideal
$I$ with the properties in the Lemma.  
\end{proof}

Our goal is now to show that there is no ideal $I$ with the properties in Lemma \ref{lem:wazzup}.
Suppose to the contrary that such an $I$ exists.  We then have a commutative diagram:
\begin{equation}
\label{eq:square}
\xymatrix {
\Gamma \ar[d]^{\rho_I} \ar@{=}[r]& \Gamma \ar[d]^{\rho_R}\\
\mathrm{GL}_q(\mathrm{W}[[t]]/I) \ar[r]&\mathrm{GL}_q(R) 
}
\end{equation}

\begin{lemma}
\label{lem:next} Suppose there is an ideal $I \subset \mathrm{W}[[t]]$ with the properties in Lemma 
$\ref{lem:wazzup}$.
Then $$I = (p^{n+1} t, pt^2,t^3, ap^n t + bt^2)$$ for some $a, b \in \mathrm{W}$ such that at least one of $a$ or $b$ is a unit.
Let $C = \mathrm{W}[[t]]/I$.  Define $S$ to be the union of $\{0\}$
with the set of Teichm\"uller lifts in $\mathrm{W} = W(k)$ of the elements of $k^*$.    
Let $\mathrm{I}_q$ be the $q\times q$ identity
matrix in $\mathrm{Mat}_q(C)$.  If $g \in K \subset \Gamma$, there is a unique
$\alpha(g) \in \mathrm{Mat}_q(S)$ and a $\xi(g) \in \mathrm{Mat}_q(\mathrm{W})$ such that in 
$\mathrm{Mat}_q(C)$ one of the following
mutually exclusive alternatives holds:
\begin{enumerate}
\item[(a)] One has that $b$ is a unit, 
$t^2=-b^{-1}ap^nt$
and $p^{n+1}t=0$ in $C$ and 
\begin{equation}
\label{eq:compute}
\rho_I(g) = \mathrm{I}_q +  t \alpha(g) + pt \xi(g).
\end{equation}  
\item[(b)] One has that $a$ is a unit in $\mathrm{W}$ and $b \in p\mathrm{W}$, 
$p^n t = 0=pt^2=t^3 $ in $C$ and  
\begin{equation}
\label{eq:compute2}
\rho_I(g) = \mathrm{I}_q +  t \alpha(g)  + t^2 \beta(g)  + pt \xi(g)
\end{equation}
for a unique $\beta(g) \in S$. 
\end{enumerate}
In both cases $(a)$ and $(b)$, one has
\begin{equation}
\label{eq:rhoRp}
\left ( \psi(g) \quad \mathrm{mod} \quad pM\right ) = \left ( \alpha(g) \quad \mathrm{mod}
\quad \mathfrak{m}_C \cdot \mathrm{Mat}_q(C) \right )
\end{equation}
where $\mathfrak{m}_C$ is the maximal ideal of $C$, $\psi:K \to M$ is the homomorphism 
in $(\ref{eq:lifter})$, and we identify the element of 
$$\mathrm{Mat}_q(C)/ (\mathfrak{m}_C \cdot \mathrm{Mat}_q(C)) = \mathrm{Mat}_q(k) = 
\mathrm{Hom}_k(V,V)$$ 
on the right side of  $(\ref{eq:rhoRp})$ with an element of the left side of $(\ref{eq:rhoRp})$ via
the identification 
$$ \mathrm{Hom}_k(V,V) = M/pM.$$
\end{lemma}

\begin{proof} By assumption $(p^n t,t^2)/I$ is isomorphic to $k$, so that $I$ contains the product
ideal $$(p^n t,t^2) \cdot (p, t) = (p^{n+1}t, pt^2, t^3)$$ in $\mathrm{W}[[t]]$.  Now $(p^n t,t^2)/(p^{n+1}t,pt^2,t^3)$
is a two-dimensional vector space over $k$ with a basis given by the classes of $p^n t$ and $t^2$.  Since $\mathrm{dim}_k((p^n t,t^2)/I) = 1$
and $(p^{n+1}t, pt^2, t^3) \subset I$,  there must be an element of $I$
of the form $ap^n t + bt^2$ in which $a, b \in \mathrm{W}$ and at least one of $a$ or $b$ is a unit.
If $g \in K$, then it follows from (\ref{eq:Rex}) that $\rho_R(g) \equiv \mathrm{I}_q$ mod 
$t \,\mathrm{Mat}_q(R)$ inside  $\mathrm{GL}_q(R)$ in (\ref{eq:lifter}).  It follows from
(\ref{eq:square}) and from $C/Ct = R/Rt = \mathrm{W}$   that 
 \begin{equation}
 \label{eq:rhog}
 \rho_I(g) \equiv \mathrm{I}_q\quad \mathrm{ mod}\quad t\,\mathrm{Mat}_q(C)  
 \quad \mathrm{for}\quad g \in K.
 \end{equation} 
 
Suppose first that $b$ is a unit. 
Then 
$$t^2=-b^{-1}ap^nt \quad\mbox{in} \quad C = \mathrm{W}[[t]]/(p^{n+1}t, pt^2, t^3, ap^n t + bt^2).$$
Hence $I=(p^{n+1}t, t^2+b^{-1}ap^nt)$, since $pt^2=-b^{-1}ap^{n+1}t\in I$ and 
$t^3=-b^{-1}ap^nt^2 \in \mathrm{W}pt^2\subset I$.
This and (\ref{eq:rhog}) lead to (\ref{eq:compute}) since 
$$C = \mathrm{W}[[t]]/I = \mathrm{W}[[t]]/(p^{n+1}t, t^2+b^{-1}ap^nt) = \mathrm{W} \oplus 
(\mathrm{W}t/\mathrm{W}p^{n+1} t).$$

Now suppose $b\in p\mathrm{W}$, so that $a$ must be a unit. Then 
$I=(p^nt,pt^2,t^3)$, since $bt^2\in \mathrm{W}pt^2$ lies in $I$, so
$(ap^n t + bt^2)-bt^2=ap^nt\in I$ and $a$ is a unit in $\mathrm{W}$. This and 
(\ref{eq:rhog}) lead to (\ref{eq:compute2}) since
$$C=\mathrm{W}[[t]]/I = \mathrm{W}[[t]]/(p^nt, pt^2,t^3) = \mathrm{W} \oplus 
(\mathrm{W}t/\mathrm{W}p^n t) \oplus (\mathrm{W}t^2/\mathrm{W}p t^2).$$

The congruence in (\ref{eq:rhoRp}) is a consequence of the fact that $\rho_I$ is a lift of $\rho_R$
(see (\ref{eq:square})) and the construction of $\rho_R$ in (\ref{eq:lifter}). 
\end{proof}

\medbreak

\noindent {\bf Completion of the proof of Theorem \ref{thm:genresult}.}
It is enough to show that there is no ideal $I \subset \mathrm{W}[[t]]$ with
the properties in Lemma \ref{lem:wazzup}.  Suppose to the contrary that such an $I$ exists.  
Let $z$ be an element of $(p^n t,t^2)$ which is not contained in the ideal $I$.  Then the
class $\overline{z}$ of $z$ in $C = \mathrm{W}[[t]]/I$ generates a $C$-ideal
$(\overline{z})$ which is the kernel of the natural surjection $C \to R$ and
which is isomorphic to $k = C/\mathfrak{m}_C$ as a $C$-module. 
We have a commutative diagram of group homomorphisms
\begin{equation}
\label{eq:bigker}
\xymatrix {
&& \Gamma \ar[d]^{\rho_I} \ar@{=}[r]&\Gamma \ar[d]^{\rho_R}\\
1\ar[r]&(1 + \overline{z}\cdot \mathrm{Mat}_q(C))^* \ar[r]&\mathrm{GL}_q(C)\ar[r] &\mathrm{GL}_q(R)\ar[r]&1
}
\end{equation}
The choice of an isomorphism $(\overline{z}) \simeq k$ of $C$-modules identifies 
the multiplicative group  \linebreak[4]
$(1 + \overline{z}\cdot \mathrm{Mat}_q(C))^*$ with $\mathrm{Hom}_k(V,V)$.

One of alternatives
(a) or (b) of Lemma \ref{lem:next} holds. 
Let $g$ be an element of $K$.   
Suppose first that $a \in p\mathrm{W}$, so that $b$ is a unit and we have alternative (a). 
 Since
 $g^{p^n}=e_K$ is the identity element of $K$, we have
\begin{equation}
\label{eq:contra}
\mathrm{I}_q = \rho_I(g)^{p^n} = (\mathrm{I}_q + t (\alpha(g) + p \xi(g)))^{p^n} = \mathrm{I}_q + p^n t  \alpha(g)
\end{equation}
since $t^2 = 0 = p^{n+1} t$ in $C$ when $a \in p\mathrm{W}$ and $b$ is a unit.  
Thus $p^n t \alpha(g) = 0$ for all $ g \in K$.  The map $p^n tC \to k = C/\mathfrak{m}_C$ defined
by $p^n tr \to (r \ \mathrm{mod} \ \mathfrak{m}_C)$ is an isomorphism because alternative (a) holds.
Thus $\alpha(g) = 0$ in $C/\mathfrak{m}_C$ for all $g \in K$.  Because of (\ref{eq:rhoRp}),
we conclude that 
$\psi(g) \in pM$ for all $g \in K$.  
Since $\psi:K=V'\to M$ is an injection and $K= V'$ and $M$ both have exponent $p^n$ as abelian
groups, this is a contradiction.

Suppose next that $a$ is a unit. Then either $b$ is a unit and we have alternative (a), which means that
$p^{n+1}t = 0$ and $t^2=-b^{-1}ap^nt$ in $C$, or $b\in p\mathrm{W}$ and we have alternative (b),
which means that $p^nt=0=pt^2=t^3$ in $C$. Suppose $h$ is another element in $K$.
If we have alternative (a), then $pt^2=0$ and we have from (\ref{eq:compute}) that 
\begin{eqnarray}
\label{eq:hgproda}
\rho_I(gh) &=& \rho_I(g)\cdot \rho_I(h)\nonumber\\
&=& (\mathrm{I}_q + t\, \alpha(g) + pt\, \xi(g) ) \cdot (\mathrm{I}_q +  t\, \alpha(h) + pt\,\xi(h))\nonumber\\
&=& \mathrm{I}_q +  t\,[\alpha(g) + \alpha(h) + p(\xi(g) + \xi(h))] + t^2\,[\alpha(g)\cdot \alpha(h)]
\nonumber\\
&=& \mathrm{I}_q +  t\,[\alpha(g) + \alpha(h) + p(\xi(g) + \xi(h))-(b^{-1}ap^n)\,\alpha(g)\cdot \alpha(h)].
\end{eqnarray}
If we have alternative (b), then $pt^2=0=t^3$ and we have from 
(\ref{eq:compute2}) that 
\begin{eqnarray}
\label{eq:hgprodb}
\rho_I(gh) &=& \rho_I(g)\cdot \rho_I(h)\nonumber\\
&=& (\mathrm{I}_q + t \,\alpha(g) + t^2 \beta(g) + pt \,\xi(g) ) \cdot (\mathrm{I}_q +  t \,\alpha(h) +
 t^2 \beta(h) + pt\,\xi(h))\nonumber\\
&=& \mathrm{I}_q +  t\,[\alpha(g) + \alpha(h) + p(\xi(g) + \xi(h))] + t^2\,[\beta(g) + \beta(h) +
\alpha(g)\cdot \alpha(h)].
\end{eqnarray}
Because $K$ is abelian, we must have $\rho_I(hg) = \rho_I(gh)$.
So subtracting (\ref{eq:hgproda}) (resp. (\ref{eq:hgprodb})) from the same formula when $g$ and $h$ are 
switched leads to the fact that 
\begin{equation}
\label{eq:switch}
\alpha(g) \cdot \alpha(h) = \alpha(h) \cdot \alpha(g) \quad \mbox{for all $g,h \in K$}
\end{equation}
for the following reason.
The map $Cp^nt \to k = C/\mathfrak{m}_C$ defined by $rp^nt \to (r \ \mathrm{mod} \ 
\mathfrak{m}_C)$ is an isomorphism when alternative (a) of Lemma \ref{lem:next} holds
and the map $C t^2 \to k = C/\mathfrak{m}_C$ defined by $rt^2 \to (r \ \mathrm{mod} \ \mathfrak{m}_C)$ is an isomorphism when alternative (b) of Lemma \ref{lem:next} holds.
Now the construction of (\ref{eq:rhoRp}) shows that $\alpha(K)$ mod $\mathfrak{m}_C\cdot
\mathrm{Mat}_q(C)$ is identified with the left $\mathbb{F}_pG$-submodule
of $M/pM = \mathrm{Hom}_k(V,V) = \mathrm{Mat}_q(k)$ which is the image of 
$\psi(K)\subset M$ in $M/pM$.  However, we assumed in Theorem \ref{thm:genresult}(iii) that 
there are $\tau, \lambda \in \psi(V')=\psi(K)$ which do not commute mod $pM$ with respect
to the product on $M  = \mathrm{Hom}_A(\hat{V},\hat{V})$ defined by the composition of maps.  This
is the product which occurs in (\ref{eq:switch}), so we have a contradiction.  
This completes the proof of  Theorem \ref{thm:genresult}.

%%%%%%%%%%%%%%%%%%%%%%%%%%%%%%%%%%%%%%%%%%%%%%%%%%%%%%%%%%%%%%%%%%%%%%%%%%%
%% Semi-direct product examples over finite fields
%%%%%%%%%%%%%%%%%%%%%%%%%%%%%%%%%%%%%%%%%%%%%%%%%%%%%%%%%%%%%%%%%%%%%%%%%%%

\section{Semi-direct product examples over finite fields}
\label{s:semdis}

Let $p$ be a prime. Suppose $q$ and $\ell$ are integers
larger than $1$ which are prime to $p$ and for which there is an injection 
$\nu:(\mathbb{Z}/q) \to (\mathbb{Z}/\ell)^*$.  We may then form the semi-direct
product $G =  (\mathbb{Z}/\ell) \semid_{\nu} (\mathbb{Z}/q)$.  We fix generators
$\tau$ and $\sigma$ for the subgroups $\mathbb{Z}/\ell$ and $\mathbb{Z}/q$ 
in this description of $G$ such that $\sigma \tau \sigma^{-1} = \tau^{\nu(\sigma)}$.

Let $\overline{\mathbb{F}}_p$ be an algebraic closure of $\mathbb{F}_p$ and
let $W(\overline{\mathbb{F}}_p)$ be the ring of infinite Witt vectors over $\overline{\mathbb{F}}_p$.
The following result is well-known so we will only sketch the proof.

\begin{lemma}
\label{lem:easyrep}
For all finite  fields $\mathbb{F}$ of characteristic $p$ the group algebra $\mathbb{F}G$ is semi-simple and
isomorphic to a direct sum of matrix algebras over finite extension fields of $\mathbb{F}$.  
Suppose $L$ is an $\overline{\mathbb{F}}_pG$-module of finite dimension over an algebraic closure $\overline{\mathbb{F}}_p$
of $\mathbb{F}_p$.
Let $B_L:G \to W(\overline{\mathbb{F}}_p)$ be the Brauer character of $L$.   
\begin{enumerate}
\item[(a)] The isomorphism class of $L$ is determined by $B_L$.
\item[(b)] Up to isomorphism, there is a unique minimal $($finite$)$ extension $k_0$ of $\mathbb{F}_p$ 
for which there is a  $k_0G$-module $L_0$ such that $L$ is isomorphic to $\overline{\mathbb{F}}_p 
\otimes_{k_0} L_0$ for some embedding of $k_0$ into $\overline{\mathbb{F}}_p$.
The field $k_0$ is isomorphic to the residue field of the ring of integers of the extension 
of $\mathbb{Q}_p$ generated by the values of $B_L$.
\item[(c)] 
The $\overline{\mathbb{F}}_pG$-module $L$ is irreducible if and only if $L_0$ from part $(b)$ is irreducible as a $k_0 G$-module. 
\item[(d)] Suppose $L_0$ in part $(b)$ is an irreducible $k_0 G$-module.  Let $L_1 = \mathrm{res}_{k_0 G}^{\mathbb{F}_pG} L_0$ 
be the 
$\mathbb{F}_pG$-module formed by the action of $\mathbb{F}_pG$ on  $L_0$.  Then $L_1$  
is irreducible as an $\mathbb{F}_pG$-module, and one has $\mathrm{End}_{\mathbb{F}_pG}(L_1) = 
\mathrm{End}_{k_0 G}(L_0) = k_0$.  
\end{enumerate}
\end{lemma}

\begin{proof}
We only show part (b).
Suppose that $k_1$ is a finite extension of $\mathbb{F}_p$ for which there
is an embedding $k_1 \to \overline{\mathbb{F}}_p$ and a 
$k_1G$-module $L_1$ such that $\overline{\mathbb{F}}_p \otimes_{k_1} L_1$
has Brauer character $B_L$.  Let $N(k_1)$ be the fraction field of the ring
$W(k_1)$ of infinite Witt vectors over $k_1$.  By \cite[\S 15.5, Prop. 43]{Serre}, there exists a 
$W(k_1)G$ module $M_1$ whose reduction mod $p$ is isomorphic to $L_1$,
and by \cite[\S 18.1(vi)]{Serre}, the character of $N(k_1) \otimes_{W(k_1)} M_1$
is $B_L$.  Thus $N(k_1)$ contains the field $N$ generated over $\mathbb{Q}_p$
by the values of $B_L$.  Hence $N$ is a finite unramified extension of $\mathbb{Q}_p$,
and is thus isomorphic to $N(k_0)$ for a unique finite extension $k_0$ of $\mathbb{F}_p$.  
To prove part (b), it will suffice to show that $L$ may be realized over $k_0$,
in the sense that there is a $k_0G$-module $L_0$ as in part (b).  By 
\cite[\S 12.2]{Serre}, there is a positive
integer $m$ such that $m\cdot B_L$ is the character of an $N(k_0)G$-module
$T$.  Then $T = N(k_0) \otimes_{W(k_0)} T_0$ for some $W(k_0)G$-module
$T_0$. The $k_0G$-module $T_1 = T_0/pT_0$ has the property that
$\overline{\mathbb{F}}_p \otimes_{k_0} T_1$ is isomorphic to a direct sum of $m$ copies
of $L$ by \cite[\S 18.2, Cor. 1]{Serre}.  Thus the class of $L$ in 
$G_0(\overline{\mathbb{F}}_pG)$ has torsion image in 
the quotient $G_0(\overline{\mathbb{F}}_pG)/ G_0(k_0G)$.  However,
the latter group is torsion free by \cite[\S 14.6]{Serre} since all finite division
rings are fields.  We thus conclude that there is a $k_0G$-module $L_0$ of the
kind required to complete the proof of part (b). 
\end{proof}

The following Lemma is a consequence of Mackey's irreducibility criterion for induced representations (see \cite[\S 7.4, Prop. 23]{Serre}).

\begin{lemma}
\label{lem:Mackey}
Let $\lambda:\langle \tau \rangle \to W(\overline{\mathbb{F}}_p)^*$ be the Brauer character
of a one-dimensional $\overline{\mathbb{F}}_p\langle\tau\rangle$-module $X_\lambda$.
Let $H = \nu(\mathbb{Z}/q) = \nu(\langle \sigma \rangle) \subset (\mathbb{Z}/\ell)^*$ be the  image
of the homomorphism $\nu:\mathbb{Z}/q = \langle \sigma \rangle \to (\mathbb{Z}/\ell)^*$ used to form $G = (\mathbb{Z}/\ell) \semid_{\nu} (\mathbb{Z}/q) = 
\langle \tau \rangle  \semid_{\nu} \langle \sigma \rangle$.  The Brauer character 
$\xi$ of the
induced $\overline{\mathbb{F}}_pG$-module $\mathrm{Ind}_{\langle \tau \rangle}^G X_\lambda$ 
has values
\begin{equation}
\label{eq:charval}
\xi (\tau^c) = \sum_{h \in H} \lambda(\tau)^{ch}\quad \mathrm{and}\quad \xi (g) = 0\quad \mathrm{if} \quad g \in G - \langle \tau \rangle
\end{equation}
for all integers $c$.   An $\overline{\mathbb{F}}_pG$-module $Y_\xi$ with Brauer character $\xi$ can be 
described as the $\overline{\mathbb{F}}_p$-vector space having a basis 
$\{w_s\}_{s \in \langle \sigma \rangle}$ on which 
$\sigma$ acts by $\sigma w_s = w_{\sigma s}$ and $\tau$ acts by 
 $\tau w_s = \overline{\lambda(\tau)}^{\nu(s^{-1})} w_s$  for all $s \in \langle \sigma \rangle$,
where $\overline{\lambda(\tau)}$ denotes the reduction of $\lambda(\tau)\in W(\overline{\mathbb{F}}_p)$
mod $p$.  The Brauer character $\xi$  is irreducible if and only if  $\lambda^h \ne \lambda$ for all 
$1 \ne h \in H$.   
\end{lemma}

\begin{cor}
\label{cor:charval}  
Suppose the Brauer character $\lambda$ of $X_\lambda$ in Lemma $\ref{lem:Mackey}$ is faithful.  
Then the Brauer character $\xi$ of  $\mathrm{Ind}_{\langle \tau \rangle}^G X_\lambda$ is faithful.  
Let $k$ be the residue field of the ring of integers of the extension of $\mathbb{Q}_p$ generated
by the values of $\xi$.
Then $k$  contains the residue field of the ring of integers of the extension of $\mathbb{Q}_p$ 
generated by the values of  the character of $\mathrm{Ind}_{\langle \tau \rangle}^G X_{\lambda^a}$ 
for all integers $a$. Finally, $k$ does not depend on the choice of the  faithful Brauer character 
$\lambda$. 
\end{cor}

\begin{proof}  
The kernel of a representation associated to $\xi=\mathrm{Ind}_{\langle \tau \rangle}^G \lambda$
must lie in $\langle \tau \rangle$, and the restriction of $\xi$ to $\langle \tau \rangle$ is a direct sum of powers of $\lambda$,
so $\xi$ is faithful because $\lambda$ is faithful.  From (\ref{eq:charval}) we see that $\xi$ and $\xi' = \mathrm{Ind}_{\langle \tau \rangle}^G \lambda^a$ take value $0$ on all $g \in G - \langle \tau \rangle$. For $g = \tau^c \in \langle \tau \rangle$ we have
$$\xi'(\tau^c) = \sum_{h \in H} \lambda^a(\tau)^{ch} = 
\sum_{h \in H} \lambda(\tau)^{ach} = \xi(\tau^{ca}).$$
Thus the values of $\xi'$ form a subset of the values of $\xi$, 
so $k$  contains the residue field of the ring of integers of the extension of $\mathbb{Q}_p$ 
generated by the values of $\xi'$. Every faithful
character of $\langle \tau \rangle $ has the form $\lambda^a$ for some $a$ 
which is relatively prime to $\ell$.
Since the map $c \to a c$ is a permutation of $\mathbb{Z}/\ell$ for such $a$, it follows 
that $k$ does not depend on the choice of $\lambda$.
\end{proof}

\begin{dfn}
\label{dfn:thegroup} 
Let $\theta: \langle \tau \rangle \to W(\overline{\mathbb{F}}_p)^*$ be a fixed faithful Brauer character
with corresponding one-dimensional $\overline{\mathbb{F}}_p\langle\tau\rangle$-module $X_\theta$, 
and suppose $a \in (\mathbb{Z}/\ell)$.  Let $k$ be the field in Corollary \ref{cor:charval} when we let 
$\lambda = \theta$.  In particular, $\mathrm{Ind}_{\langle \tau \rangle}^G \theta: \langle \tau \rangle \to 
W(k)$ when $W(k)$ is the ring of 
infinite Witt vectors over $k$. Let $V(\theta^a)$ be a $kG$-module with Brauer character 
$B_{V(\theta^a)} = \mathrm{Ind}_{\langle \tau \rangle}^G \theta^a:G\to W(k)$.  
(Note that $\theta^a$ need not be faithful, so that 
the residue field of the ring of integers of the extension of $\mathbb{Q}_p$ 
generated by the values of $B_{V(\theta^a)}$ need not be $k$.)
Fix an integer $n \ge 1$, and let $A = W(k)/p^n W(k)$. Let $\hat{V}(\theta^a)$ be an $AG$-module
which is finitely generated and free as an $A$-module such that $k \otimes_A \hat{V}(\theta^a)$ is 
isomorphic to $V(\theta^a)$.  
Let $K=\hat{V}(\theta^a)$ as abelian $p$-groups and let $\delta:G \to \mathrm{Aut}(K)$ be the 
group homomorphism given by the $G$-action on the left $AG$-module  $\hat{V}(\theta^a)$.
Define $\Gamma$ to be the semi-direct product $K \semid_{\delta} \,G$.
If $L$ is a $G$-module, we let $\tilde L$ be the $\Gamma$-module which is the
inflation of $L$ via the natural surjection $\pi:\Gamma \to G$.  Let $V = V(\theta)$, and let $\hat{V} = \hat{V}(\theta)$.
Define 
$M = \mathrm{Hom}_A(\hat{V},\hat{V})$ and $M_0 = k\otimes_A M = \mathrm{Hom}_k(V,V)$ with the usual $G$-actions 
coming from the action of $G$
on $\hat{V}$ and $V$.  We will sometimes fix an $A$-basis of $\hat{V}$ and thus an isomorphism of 
$M$ with the matrix algebra $\mathrm{Mat}_q(A)$ and an isomorphism of $M_0$ with $\mathrm{Mat}_q(k)$.  
\end{dfn}

We will now make the following assumptions:

\begin{hypo}
\label{hyp:Hneeded}
Assume the notation of Definition $\ref{dfn:thegroup}$.  
Suppose $0 \ne a \in (\mathbb{Z}/\ell)$ and that $H = \nu(\mathbb{Z}/q) = \nu(\langle \sigma \rangle) \subset (\mathbb{Z}/\ell)^*$ as in Lemma $\ref{lem:Mackey}$.   Suppose further that
\begin{enumerate}
\item[(a)] $(h-1)a \ne 0$ in $\mathbb{Z}/\ell$ for $1 \ne h \in H$, and that
\item[(b)] there is
exactly one ordered pair $(h_2, h_3)$ of elements 
$h_2, h_3 \in H$ such that $h_3 - h_2 = a$ in $\mathbb{Z}/\ell$, and that
\item[(c)] the residue field of the ring of integers of the extension of $\mathbb{Q}_p$ generated
by the values of the Brauer character $B_{V(\theta^a)}$ is equal to $k$.
\end{enumerate}
\end{hypo}

\begin{thm}
\label{thm:main} 
Suppose Hypothesis $\ref{hyp:Hneeded}$ holds.  Then all the hypotheses of Theorem 
$\ref{thm:genresult}$ hold when
we let $V = V(\theta)$ and we let $V' = \mathrm{res}_{AG}^{(\mathbb{Z}/p^n)G} \hat{V}(\theta^a)$ be the 
$(\mathbb{Z}/p^n)G$-module which results from $\hat{V}(\theta^a)$ by
restricting operators from $AG$ to $(\mathbb{Z}/p^n)G$. 
In particular, $\mathrm{End}_{kG}(V)=k$. It follows that the universal
deformation ring  $R_{W(k)}(\Gamma,\tilde {V})$ of $\tilde {V}$  is isomorphic to $W(k)[[t]]/(p^n t,t^2)$.
\end{thm}

The proof of this result will occupy the rest of this section.  Let us first show that Theorem \ref{thm:main} implies Theorem \ref{thm:vaguemain}
by showing that the hypotheses of Corollary \ref{cor:reduce} may always be satisfied.
If $p = 2$, this is shown by Example \ref{eq:easyex}.  

Suppose $p = 3$.   Let  $q=2$, $\ell = 8$, $H = \{1, 3\} \subset (\mathbb{Z}/8)^*$ and $a= 2$.  
The unique element $h$ of $H$ which  is not $1$ is $h  = 3$, and $(3 -1) a=4 \ne 0$ in 
$\mathbb{Z}/\ell  = \mathbb{Z}/8$.
The unique ordered pair $(h_2,h_3)$ of elements of $H$
for which $h_2 - h_3 = a$ in $\mathbb{Z}/\ell$ is $(h_2,h_3) = (3, 1)$.  Thus Hypothesis 
\ref{hyp:Hneeded} holds.
Note that (\ref{eq:charval}) and Corollary \ref{cor:charval} show that 
the residue field $k$ of the ring of integers of the extension of $\mathbb{Q}_3$ 
generated by the values of the Brauer character $\mathrm{Ind}_{\langle \tau \rangle}^G \theta$
is the extension of $\mathbb{F}_3$ generated by $\zeta + \zeta^3$ as $\zeta$ ranges over all roots of unity of order dividing  $\ell=8$  in $\overline{\mathbb{F}}_3$.
Since the absolute Frobenius automorphism  $\alpha \to \alpha^3$ fixes $\zeta+\zeta^3$,
we conclude that $k = \mathbb{F}_3$.   
Thus Theorem \ref{thm:main} and Corollary \ref{cor:reduce} 
show Theorem \ref{thm:vaguemain} for $p = 3$.    

Suppose now that $p > 3$.  Define $q = 2$ and $\ell = 3$.  Then $H = \{\pm 1\} = (\mathbb{Z}/\ell)^*$. 
Let $a = 1$. We readily check that Hypothesis \ref{hyp:Hneeded} holds.  
The field $k$
is the extension of $\mathbb{F}_p$ generated by $\zeta + \zeta^{-1}$ 
as $\zeta$ ranges over all roots of unity of order  dividing $\ell = 3$ in 
$\overline{\mathbb{F}}_p$.  Thus  $k = \mathbb{F}_p$, so Theorem \ref{thm:main} 
and Corollary \ref{cor:reduce} show Theorem \ref{thm:vaguemain} for $p > 3$.

We now come back to the proof of Theorem \ref{thm:main},
which is a consequence of the following result.
Note that $V(\theta)$ is irreducible by Lemma \ref{lem:Mackey} and $\mathrm{End}_{kG}(V(\theta))=k$
by Lemma \ref{lem:easyrep}.

\begin{lemma}
\label{lem:multone}
Assume the notation of Definition $\ref{dfn:thegroup}$, so that $V = V(\theta)$, and suppose 
$a \in \mathbb{Z}/\ell$.
\begin{enumerate}
\item[(i)]  The $kG$-module $V(\theta^a)$ is irreducible if and only if $a$  satisfies condition $(a)$ 
of Hypothesis $\ref{hyp:Hneeded}$,
and in this case $V(\theta^a)$ is absolutely irreducible in the sense that 
$\overline{\mathbb{F}}_p \otimes_k V(\theta^a) $ is irreducible as an $\overline{\mathbb{F}}_pG$-module.  
\item[(ii)]  Suppose $V(\theta^a)$ is an irreducible $kG$-module.  The multiplicity of $V(\theta^a)$ in 
$M_0 = \mathrm{Hom}_k(V,V)$ equals $1$ if and only if
$a$ satisfies condition $(b)$ of Hypothesis $\ref{hyp:Hneeded}$.
\end{enumerate}
Suppose now that all three conditions $(a)$, $(b)$ and $(c)$ of Hypothesis $\ref{hyp:Hneeded}$ hold.  Then:
\begin{enumerate}
\item[(iii)] 
If $V'=\mathrm{res}_{AG}^{(\mathbb{Z}/p^n)G} \hat{V}(\theta^a)$, then
$\mathrm{Hom}_{(\mathbb{Z}/p^n)G}(V',M)$ is a free rank one $A$-module with
respect to the $A$-module structure coming from the multiplication action of $A$ on $M$.
There is an injective $AG$-module homomorphism $\psi:\hat{V}(\theta^a) \to M$ which is unique up 
to multiplication by an element of $A^*$ and which defines
an injective homomorphism in $\mathrm{Hom}_{(\mathbb{Z}/p^n)G}(V',M)$.
\item[(iv)] Relative to the ring structure for $M = \mathrm{Hom}_A(\hat{V},\hat{V})$ coming from the 
composition of homomorphisms, there are elements of $\psi(\hat{V}(\theta^a))$ which
do not commute mod $pM$.
\end{enumerate}  
\end{lemma}

\begin{proof}  By Lemma \ref{lem:Mackey}, $\overline{\mathbb{F}}_p \otimes_k V(\theta^a)$ is  
irreducible if and only if $\theta^{ah} \ne \theta^a$ for all $1 \ne h \in H$.  Since $\theta$ is a faithful 
Brauer character of $\mathbb{Z}/\ell$,
this is the case if and only if $(h-1)a \ne 0$ in $\mathbb{Z}/\ell$ for $1 \ne h \in H$, which is condition (a)
of Hypothesis \ref{hyp:Hneeded}. Clearly $V(\theta^a)$ is irreducible as a $kG$-module if it is absolutely
irreducible.  Conversely, suppose $V(\theta^a)$ is an irreducible $kG$-module.  By Corollary
\ref{cor:charval}, $k$ contains the 
residue field $k_0$ of the ring of integers of the extension of $\mathbb{Q}_p$ 
generated by the values of the Brauer character $B_{V(\theta^a)}:G\to W(k)$. 
By Lemma \ref{lem:easyrep}, $\overline{\mathbb{F}}_p \otimes_k V(\theta^a) $ is isomorphic as an 
$\overline{\mathbb{F}}_pG$-module
to $\overline{\mathbb{F}}_p \otimes_{k_0} L_0$ for some $k_0G$-module $L_0$.  Since $V(\theta^a)$ 
has the same Brauer character 
as $k \otimes_{k_0} L_0$, it follows from the semi-simplicity of $kG$ that
$V(\theta^a)$ is isomorphic to $k \otimes_{k_0} L_0$.  Hence $L_0$ must be irreducible
as a $k_0G$-module since $V(\theta^a)$ is irreducible as a $kG$-module.  Now
$\overline{\mathbb{F}}_p \otimes_k V(\theta^a)  \simeq \overline{\mathbb{F}}_p \otimes_{k_0} L_0$ 
is irreducible as an $\overline{\mathbb{F}}_pG$-module by Lemma \ref{lem:Mackey}(c).

To compute the multiplicity of $V(\theta^a) $ in $M_0$ we will
use Brauer characters.
By (\ref{eq:charval}), we have for all integers $b$ and $c$ that 
\begin{equation}
\label{eq:ohyeah}
B_{V(\theta^b)}(\tau^c) = \sum_{h \in H} \zeta^{bch}\quad \mathrm{and}\quad 
B_{V(\theta^b)}(g) = 0\quad \mathrm{if} \quad g \in G - \langle \tau \rangle
\end{equation}
where $\zeta=\theta(\tau)$ is a primitive $\ell^{\mathrm{th}}$ root of unity in $W(\overline{\mathbb{F}}_p)$.
Let $J^* = \mathrm{Hom}_k(J,k)$ be the contragredient of
a finitely generated $kG$-module $J$.   The Brauer character $B_{J^*}$ satisfies 
$B_{J^*}(g) = B_{J}(g^{-1})$ for all $g\in G$.  
There is a $kG$-module isomorphism 
$$M_0 = \mathrm{Hom}_k(V,V) \cong V^* \otimes_k V .$$
The Brauer character $B_{M_0}$ of $M_0$ is thus  $B_{V}\cdot B_{V^*}$.

We now calculate the multiplicity of the irreducible representation $V(\theta^a)$
in $M_0$ using the standard inner product on Brauer characters.  This gives:
\begin{eqnarray}
\label{eq:Mokey}
\mathrm{mult}(V(\theta^a),M_0) & = & \frac{1}{\#G} \sum_{g \in G} B_{V(\theta^a)}(g) \cdot 
B_{M_0}(g^{-1})\nonumber\\
&=& \frac{1}{q\ell} \sum_{c = 0}^{\ell-1}  B_{V(\theta^a)}(\tau^c) \cdot B_{M_0}(\tau^{-c})\nonumber\\
&=&\frac{1}{q\ell} \sum_{c = 0}^{\ell -1} \left(\sum_{h_1 \in H} \zeta^{ach_1}\right ) \cdot \left (\sum_{h_2 \in H} \zeta^{ch_2}\right ) \cdot \left (\sum_{h_3 \in H} \zeta^{-ch_3} \right )\nonumber\\
&=&\frac{1}{q\ell} \sum_{c = 0}^{\ell -1} \sum_{h_1, h_2, h_3 \in H} \zeta^{c(ah_1+ h_2 - h_3)}
\end{eqnarray}
Since $\zeta$ is a primitive $\ell^{\mathrm{th}}$ root of unity in $W(\overline{\mathbb{F}}_p)$, we have 
$\sum_{c = 0}^{\ell -1} \zeta^{cd} = 0$ if $d \not \equiv 0$ mod $\ell$
and $\sum_{c = 0}^{\ell -1} \zeta^{cd} = \ell$ if $d \equiv 0$ mod $\ell$.  So if we let $S$ be the set of all ordered triples $(h_1, h_2, h_3)$
of elements of $H$ such that $ah_1 + h_2 - h_3 = 0$ in $\mathbb{Z}/\ell$, we see from 
(\ref{eq:Mokey}) that 
\begin{equation}
\label{eq:mo2}
\mathrm{mult}(V(\theta^a),M_0) = \frac{\#S}{q}.
\end{equation}
We now observe that $h \in H$ acts on $S$ by sending $(h_1, h_2, h_3)$ to $(hh_1, hh_2, hh_3)$.  Since $H$ is a group
of order $q$, we see that $\mathrm{mult}(V(\theta^a),M_0) = 1$ in (\ref{eq:mo2}) if and only if there is a unique triple of the
form $(1,h_2,h_3)$ in $S$.  This is equivalent to the statement that there is a unique ordered pair $(h_2,h_3)$ of elements of $H$
such that $a \equiv h_3 - h_2$ mod $\ell$. 

Suppose now that conditions (a), (b) and (c) of Hypothesis \ref{hyp:Hneeded} are satisfied. Then 
$V(\theta^a)$ is absolutely irreducible by part (i) and it has multiplicity 1 in $M_0$ by part (ii).
Since $kG$ is semi-simple, this implies that $\mathrm{Hom}_{kG}(V(\theta^a),M_0) \cong k$
and that there is an injective $kG$-module homorphism $\psi_0:V(\theta^a) \to M_0$.
Let $V'=\mathrm{res}_{AG}^{(\mathbb{Z}/p^n)G} \hat{V}(\theta^a)$, which means that
$V'/pV'=\mathrm{res}_{kG}^{\mathbb{F}_pG} V(\theta^a)$.
By condition (c) of Hypothesis \ref{hyp:Hneeded}, the residue field of the ring of integers of the extension of 
$\mathbb{Q}_p$ generated by the values of the Brauer character $B_{V(\theta^a)}:G\to W(k)$ is equal to 
$k$. By Lemma \ref{lem:easyrep}, this implies that $V'/pV'$ is a simple
$\mathbb{F}_pG$-module with $\mathrm{End}_{\mathbb{F}_pG}(V'/pV')\cong k$. Since
$V(\theta^a)$ has multiplicity 1 in $M_0$, it follows that $\mathrm{Hom}_{\mathbb{F}_pG}(V'/pV',M_0)
\cong \mathrm{End}_{\mathbb{F}_pG}(V'/pV')\cong k$.

Recall that $A = W(k)/p^nW(k)$ and that $\hat{V}(\theta^a)$ and $M$ are
projective $AG$-modules such that $\hat{V}(\theta^a)/p\hat{V}(\theta^a)$ is isomorphic to $V(\theta^a)$
and $M/pM$ is isomorphic to $M_0$ as $kG$-modules.  It follows that 
$\mathrm{Hom}_{AG}(\hat{V}(\theta^a),M)$ is a projective $A$-module $T$ such that
$T/pT = \mathrm{Hom}_{kG}(V(\theta^a),M_0) \cong k$.  Thus  
$\mathrm{Hom}_{AG}(\hat{V}(\theta^a),M) \cong A$,
which implies that there is an injective $AG$-module homomorphism $\psi:\hat{V}(\theta^a)\to M$ 
which is unique up to multiplication by an element of $A^*$.
Since $\mathrm{Hom}_{AG}(\hat{V}(\theta^a),M)\cong A$ is an $A$-submodule of 
$\mathrm{Hom}_{(\mathbb{Z}/p^n)G}(V',M)$ and since the reductions mod $p$
of both these $A$-modules are isomorphic to $k$, it follows that 
$\mathrm{Hom}_{(\mathbb{Z}/p^n)G}(V',M)\cong A$.

We still need to show that $\psi(\hat{V}(\theta^a))$ is a non-commutative subset of 
$M = \mathrm{Hom}_A(\hat{V},\hat{V})$ mod $pM$ with respect to the multiplication
coming from the composition of homomorphisms.  We can identify $M/pM$ with $M_0$ and the image of $\psi(\hat{V}(\theta^a))$
in $M/pM = M_0$ with $\psi_0(V(\theta^a))$.  So we must show $\psi_0(V(\theta^a))$ is a non-commutative subset of
$M_0 = \mathrm{Hom}_{kG}(V,V)$.   

One has 
\begin{equation}
\label{eq:dimone}
\overline{\mathbb{F}}_p \otimes_k \mathrm{Hom}_{kG}(V(\theta^a),M_0) = \mathrm{Hom}_{\overline{\mathbb{F}}_pG}(
\overline{\mathbb{F}}_p \otimes_k V(\theta^a), \overline{\mathbb{F}}_p \otimes_k M_0)
\end{equation}
where
$$\overline{\mathbb{F}}_p \otimes_k M_0 = \mathrm{Hom}_{\overline{\mathbb{F}}_pG}(
\overline{\mathbb{F}}_p \otimes_k V,\overline{\mathbb{F}}_p \otimes_k V).$$
{}From Lemma \ref{lem:Mackey} we have isomorphisms
\begin{equation}
\label{eq:Vades}
\overline{\mathbb{F}}_p \otimes_k V(\theta^a)  = \bigoplus_{s \in \langle \sigma \rangle} 
\overline{\mathbb{F}}_p w_s
\end{equation}
and
\begin{equation}
\label{eq:Vdes}
\mathrm{Hom}_{\overline{\mathbb{F}}_p}(\overline{\mathbb{F}}_p \otimes_k V,\overline{\mathbb{F}}_p \otimes_k V) =
\bigoplus_{s',s'' \in \langle \sigma \rangle} \mathrm{Hom}_{\overline{\mathbb{F}}_p}(\overline{\mathbb{F}}_p x_{s'},\overline{\mathbb{F}}_p x_{s''})
\end{equation}
where  for $s \in \langle \sigma \rangle$ we have  
\begin{equation}
\label{eq:tauform}
\tau w_s =    \overline{\theta^a(\tau)}^{\nu(s^{-1})} w_s\quad \mathrm{and} \quad 
\tau x_s = \overline{\theta(\tau)}^{\nu(s^{-1})} x_s
\end{equation}
and
\begin{equation}
\label{eq:sigmaform}
\sigma w_s = w_{\sigma s} \quad \mathrm{and}\quad \sigma x_s = x_{\sigma s}.
\end{equation}

We now exhibit a non-zero $\overline{\mathbb{F}}_pG$-module homomorphism $f$ from
(\ref{eq:Vades}) to (\ref{eq:Vdes}).  Recall that there is a unique ordered pair $(h_2, h_3)$ of elements 
$h_2, h_3 \in H$ such that $h_3 - h_2 = a$ in $\mathbb{Z}/\ell$.  Suppose $s \in \langle \sigma \rangle$.
Let $s_2, s_3 \in \langle \sigma \rangle$
be the unique elements of $\langle \sigma \rangle$ such that $\nu(s_2) = \nu(s) h_2^{-1}$ and $\nu(s_3) = \nu(s) h_3^{-1}$.
Define $f(w_s)$ to be the unique element of (\ref{eq:Vdes}) which sends $x_{s_2}$ to $x_{s_3}$ and which sends
$x_{s'}$ to $0$ if $s_2 \ne s' \in \langle \sigma \rangle$.  We can then extend $f$ to a unique
$\overline{\mathbb{F}}_p$-linear map from (\ref{eq:Vades}) to (\ref{eq:Vdes}).  A  straightforward
calculation using (\ref{eq:tauform}) and (\ref{eq:sigmaform}) shows that $f$ is $G$-equivariant.  

Consider now the image of $f$. Let $e$ be the identity element of $\langle \sigma \rangle$.  Then
$f(w_e)(x_{s_2}) = x_{s_3}$ if $\nu(s_2) = h_2^{-1}$ and $\nu(s_3) = h_3^{-1}$.  There is
a unique $s \in \langle \sigma \rangle$
such that $\nu(s) = h_2 h_3^{-1}$.  Then 
\begin{equation}
\label{eq:fws}
f(w_s)(x_{s'_2}) = x_{s'_3}\quad \mathrm{when} \quad \nu(s'_2) = \nu(s) h_2^{-1}\quad \mathrm{and}\quad 
\nu(s'_3) = \nu(s) h_3^{-1}
\end{equation}
and
\begin{equation}
\label{eq:fws2}
f(w_s)(x_{s'}) = 0 \quad \mathrm{if} \quad s' \ne s'_2.
\end{equation}
  Here $\nu(s'_2) = \nu(s) h_2^{-1} = h_2 h_3^{-1} h_2^{-1} = h_3^{-1}$.
Thus $s'_2 = s_3 \ne s_2$ since $h_2$ and $h_3$ are distinct because $h_3 - h_2 =a\ne 0$ in $\mathbb{Z}/\ell$.  So
(\ref{eq:fws}) and (\ref{eq:fws2}) give 
\begin{equation}
\label{eq:fw}
(f(w_s) \circ f(w_e))(x_{s_2}) = f(w_s) (f(w_e)(x_{s_2})) = f(w_s)(x_{s_3}) = x_{s'_3} \quad \mathrm{and}\quad f(w_s)(x_{s_2}) = 0
\end{equation}
where $s'_3$ is determined by $\nu(s'_3) = \nu(s) h_3^{-1} = h_2 h_3^{-2}$.    The second statement in (\ref{eq:fw}) gives 
\begin{equation}
\label{eq:fw2}
(f(w_e) \circ f(w_s))(x_{s_2}) = f(w_e)(f(w_s)(x_{x_2})) = f(w_e)(0) = 0.
\end{equation}
Comparing (\ref{eq:fw}) and (\ref{eq:fw2}) shows that $f(w_s) \circ f(w_e) \ne f(w_e) \circ f(w_s)$,
so that $f(\overline{\mathbb{F}}_p\otimes_k V(\theta^a))$ is a non-commutative
subset of $\overline{\mathbb{F}}_p \otimes_k M_0 = \overline{\mathbb{F}}_p \otimes_k \mathrm{Hom}_k(V,V)$.

Since $V(\theta^a)$ is irreducible as a $kG$-module and occurs with multiplicity 
$1$ in $M_0$ and since $\mathrm{End}_{kG}(V(\theta^a)) = k$ by Lemma \ref{lem:easyrep},
we have $\mathrm{dim}_k(\mathrm{Hom}_{kG}(V(\theta^a),M_0)) = 1$.
Therefore (\ref{eq:dimone}) implies
$$f = \beta \otimes \psi_0$$
for some $0 \ne \beta \in \overline{\mathbb{F}}_p^*$ since $f$ and  $\psi_0:V(\theta^a) \to M_0$ are 
non-zero.  Hence 
$$f(\overline{\mathbb{F}}_p \otimes_k V(\theta^a)) = \overline{\mathbb{F}}_p  \otimes_k 
\psi_0(V(\theta^a)) \subset \overline{\mathbb{F}}_p \otimes_k M_0.$$
  Thus 
$\psi_0(V(\theta^a))\cong 1 \otimes \psi_0(V(\theta^a))$ must be non-commutative because $f(\overline{\mathbb{F}}_p \otimes_k V(\theta^a))$ is, and this completes the proof.
\end{proof}


\begin{thebibliography}{99}

\bibitem{FB1} F.~M.~Bleher, Universal deformation rings and dihedral defect groups. 
	Trans. Amer. Math. Soc. 361 (2009), 3661--3705. 

\bibitem{FB2} F.~M.~Bleher, Universal deformation rings and generalized quaternion
	defect groups. Submitted for publication in 2009. {\tt ArXiv:0909.1031}

\bibitem{BC} F.~M.~Bleher and T.~Chinburg, Universal deformation rings
	and cyclic blocks. Math. Ann. 318 (2000), 805--836.

\bibitem{lcicomptes} F.~M.~Bleher and T.~Chinburg, Universal deformation rings need not
	be complete intersections. C. R. Math. Acad. Sci. Paris 342 (2006), 229--232.

\bibitem{lcann} F.~M.~Bleher and T.~Chinburg, Universal deformation rings need not be complete 			intersections. Math. Ann. 337 (2007), 739--767.

\bibitem{JB} J.~Byszewski, A universal deformation ring which is not a complete intersection ring. 
	 C. R. Math. Acad. Sci. Paris  343  (2006),  565--568. 

\bibitem{flach} T.~Chinburg, Can deformation rings of group representations not be 
	local complete intersections? In: Problems from the Workshop on Automorphisms of Curves. 
	Edited by Gunther Cornelissen and Frans Oort, with contributions by I. Bouw, T. Chinburg, 
	Cornelissen, C. Gasbarri, D. Glass, C. Lehr, M. Matignon, Oort, R. Pries and S. Wewers.
	Rend. Sem. Mat. Univ. Padova 113 (2005), 129--177.
	
\bibitem{Darmon} H.~Darmon, F.~Diamond and R.~Taylor, Fermat's Last Theorem. In : R. 
	Bott, A. Jaffe and S. T. Yau (eds), Current developments in mathematics, 1995, International 
	Press, Cambridge, MA., 1995, pp. 1--107.
	
\bibitem{desmitlenstra} B.~de Smit and H.~W.~Lenstra, Explicit construction of universal deformation 
	rings. In: G. Cornell, J. H. Silverman and G. Stevens (eds), Modular Forms and Fermat's Last 
	Theorem (Boston, MA, 1995), Springer-Verlag, Berlin-Heidelberg-New York, 1997, pp. 313--326.
	
\bibitem{ega44} A.~Grothendieck, \'{E}l\'{e}ments de g\'{e}om\'{e}trie alg\'{e}brique, Chapitre IV, Quatri\'{e}me Partie. Publ. Math. IHES 32 (1967), 5--361.

\bibitem{matsumura} H.~Matsumura, Commutative Ring Theory. Cambridge Studies
in Advanced Mathematics, Vol. 8, Cambridge University Press, Cambridge, 1989.

\bibitem{maz1} B.~Mazur, Deforming Galois representations. In: Galois groups over $\mathbb{Q}$ 
	(Berkeley, CA, 1987), Springer-Verlag, Berlin-Heidelberg-New York, 1989, pp. 385--437.

\bibitem{maz2} B.~Mazur,  An introduction to the deformation theory of Galois representations.
	 In: G. Cornell, J. H. Silverman and G. Stevens (eds), Modular Forms and Fermat's Last 
	Theorem (Boston, MA, 1995), Springer-Verlag, Berlin-Heidelberg-New York, 1997, pp. 243--311.

\bibitem{Sch} M.~Schlessinger, Functors of Artin Rings. Trans. of the AMS
		130 (1968), 208--222.

\bibitem{Serre} J.~P.~Serre, Corps Locaux. Hermann, Paris, 1968.

\bibitem{Wiles} A.~Wiles, Modular elliptic curves and Fermat's last theorem.
        Ann. of Math. 141 (1995), 443--551.





\end{thebibliography}
\end{document}